\documentclass[12pt]{article}
\usepackage{amsmath,amsfonts,amssymb,amsthm}

\def\noi{\noindent}

\def\pf{\noi{\bf Proof.\ \,}}
\def\eop{{$\square$}}

\def\labtt#1{\label {#1}  }
\def\labttr#1{\label {#1}  \rm }

\def\a{\alpha}

\def\l{\lambda}

\def\CC{{\mathbb C}}
\def\FF{{\mathbb F}}
\def\QQ{{\mathbb Q}}
\def\RR{{\mathbb R}}

\def\ZZ{{\mathbb Z}}

\def\la{\langle}
\def\ra{\rangle}
\def\<{\langle}
\def\>{\rangle}

\def\l{{\lambda}}

\def\half{{1 \over 2}}
\def\fourth{{1 \over 4}}
\def\eighth{{1 \over 8}}
\def\sixteenth{{1 \over 16}}

\def\dual#1{#1^*}        

\def\kron#1#2{\delta_{#1#2}}  

\def\hame#1{{\cal H}_{#1}^e}

\def\leh{L_{E_8}}
\def\rtleh{\sqrt 2 L_{E_8}}

\def\mv#1{MinVec(#1)}

\def\weh{W_{E_8}}

\def\exsp#1#2{2^{1+2#1}_#2} 
\def\explus#1{2^{1+2#1}_+}

\def\ratholoex#1{2^{1+2#1}_+\Omega^+(2#1,2)}

\def\dg#1{{\cal D}({#1})}  

\def\vac{\hbox{\bf 1}} 

\def\allones#1{\hbox{\bf 1}^{#1}}

\def\bw#1{BW_{{2^{#1}}}}
 
\def\tbw{type BW}

\def\gd#1{G_{2^{#1}}} 
\def\rd#1{R_{2^{#1}}}

\begin{document}

\newtheorem{thm}{Theorem}[section]
\newtheorem{prop}[thm]{Proposition}
\newtheorem{lem}[thm]{Lemma}
\newtheorem{rem}[thm]{Remark}
\newtheorem{coro}[thm]{Corollary}
\newtheorem{conj}[thm]{Conjecture}
\newtheorem{de}[thm]{Definition}
\newtheorem{hyp}[thm]{Hypothesis}

\newtheorem{nota}[thm]{Notation}
\newtheorem{ex}[thm]{Example}
\newtheorem{proc}[thm]{Procedure}

\centerline{\Large \bf Pieces of $2^d$:   }

\centerline{\Large \bf Existence and uniqueness for }

\centerline{\Large \bf 
Barnes-Wall and    Ypsilanti lattices.  }

\centerline{(today is \today \ in Ann
Arbor.)  }
\begin{center}
{\Large    }

\vspace{10mm}
Robert L.~Griess Jr.
\\[0pt]
Department of Mathematics\\[0pt] University of Michigan\\[0pt]
Ann Arbor, MI 48109  \\[0pt]
\vskip 1cm 
Dedicated to Donald G. Higman.  
\end{center}

\vfill \eject 

\begin{abstract}
We give a new existence proof for the rank $2^d$ even lattices
usually called the Barnes-Wall lattices, and establish new results on
uniqueness, structure and transitivity of the automorphism group on
certain kinds of sublattices.  Our  
proofs are relatively free of calculations, matrix work and
counting, due to the uniqueness viewpoint.  We deduce the labeling of
coordinates on which earlier constructions depend. 

Extending these ideas, we construct in dimensions $2^d$, for $d>>0$,
{\it the Ypsilanti lattices}, which are families of indecomposable
even unimodular lattices which resemble the Barnes-Wall lattices. 
The number $\Upsilon (2^d)$ of isometry types here is large: 
$log_2 (\Upsilon (2^d))$  has dominant term at least ${r \over 4}d\,
2^{2d}$, for any $r \in [0,\half )$.  
Our lattices  may be the first explicitly given
families  whose sizes  are asymptotically comparable to the
Siegel mass formula estimate ($log_2(mass(n))$ has dominant term 
${1\over 4} log_2(n) n^2$).

This work continues our general uniqueness program for lattices, 
begun in
Pieces of Eight \cite{POE}.  See also our new uniqueness proof for the
$E_8$-lattice \cite{GrE8}.
\end{abstract}

\vfill \eject

\tableofcontents


\section{\bf Notation and terminology }

\bigbreak
\halign{#\hfil&\quad#\hfil\cr

annihilator, self annihilating  & Section 4\cr

$A_{ij}$ and other  diagonal notation  &\ref{diagonalnotation}\cr 

admissible & \ref{admissible} \cr

ancestors and generations, ancestral & \ref{ancestral},
\ref{ancestorlookalike}
\cr

$\bw d$, the Barnes-Wall lattice in dimension $2^d$ & \ref {bwt} \cr

lattice of BW-type & \ref{bwt} \cr

 $BRW^0(2^d,\pm )$ & Bolt, Room and Wall group, \ref{brw}   \cr

classification &\ref{bwtclassification} \cr 

coelementary abelian subgroup, $p$-coelementary abelian  & a subgroup
$B\le A$ so
\cr 

& 
$A/B$ is 
$p$-elementary abelian 
\cr 

$D$, a lower dihedral group  & \ref{nextbw} \cr 

defect of an involution & \ref{defectineq} \cr

density, commutator density & \ref{commutatordense},\ref{commdefs} \cr

$\dg L$, discriminant group of a lattice $L$ & $\dg L = \dual L / L$\cr

 determinant of a lattice, $L$    & $|\dg L |$ 
\cr 

$d$-invariant & \ref{dinvariant} \cr 
duality level &   \ref{dualitylevel}   \cr

double basis & \ref{lbc} \cr 

DT,  DTL & \ref{dtl}  \cr 

eigenlattice, total eigenlattice, $Tel$ & \ref{eigenlattice} \cr 

$f, f_i, f_{12}$;  various fourvolutions  &  \ref{nextbw}  \cr 

$F_i$ & \ref{nextbw} \cr

fourvolution & \ref{fourvolution} \cr

frame, 
plain frame PF, 
sultry frame SF   & \ref{plainframe}, \ref{sf} \cr 

$\gd d$ & \ref{rdgd} \cr

Hamming codes & \ref{hcsc},\ref{xhcsc} \cr

$I(d,p,q)$ & \ref{ioexp}\cr

labeling & \ref{labeling} \cr 

lower & \ref{lowerstuff}  \cr

mass formula, 
$mass(n)$  & \ref{massnota} \cr 

$L$, an integral lattice of rank $n$ & Section 5 \cr

$L_i,  L_i[k]$ & \ref{nextbw} \cr 

$\dual L$, the dual of the lattice $L$ & Section 5\cr 

$\varepsilon_A$, sign changes, monomial group & \ref{sumsandsigns} \cr 

$M_i, M_i[k]$ & \ref{nextbw} \cr 

minimal vectors, $\mv L$, $\mu (L)$  & \ref{minvec} \cr 

nextbw  & \ref{nextbw}  \cr 

power set, even sets & Section 4\cr 

$Q_i$ & \ref{nextbw} \cr 

$r$-modular & \ref{rmodlatt} \cr 

$R$ & \ref{rdgd}, \ref{nextbw} \cr 

$\rd d$ & \ref{rdgd} \cr 

$\rd i$ & \ref{nextbw} \cr 

$Scalar(G,M)$, the scalar subgroup & \ref{scalarsg} \cr  

SSD, semiselfdual,  
RSSD, relatively semiselfdual   & \ref{ssd} \cr 

sultry frame & \ref{sf} \cr 

sultry transformation, twist  & \ref{twist} \cr

sBW, ssBW   & \ref{ssbw} \cr

$t_i, t_{ij}, t_{ij'}$ & \ref{nextbw} \cr

upper & \ref{lowerstuff} \cr 

Condition $X(2^d)$ & \ref{condx} \cr

$\frak X$ & \ref{frakx}\cr

$\frak Y$ & \ref{fraky} \cr 

Ypsilanti lattices, cousins, etc.   &Section 14
\cr

zop2 , 
zoop2  & \ref{zoop2} \cr

$\varepsilon_A$, sign changes, monomial group & \ref{sumsandsigns} \cr 

$\psi_i$ & \ref{nextbw} \cr 

$\Omega$, index set (often identified with an affine space
$\FF_2^d$)&  \cr 

$\Omega$, universe, $\vac$, the ``all ones vector'' $(1,1,\dots,1)$ in
$\FF_2^\Omega$ & Section 4
\cr

\cr }

{\bf Conventions.  }  Our groups and most endomorphisms act on the
right, often with exponential notation.  
Group  theory notation is
mostly consistent with \cite{Gor, Hup, G12}. 
The commutator of $x$ and $y$ means $[x,y]=x^{-1}y^{-1}xy$ and
the  conjugate of of
$x$ by
$y$ means
$x^y:=y^{-1}xy=x[x,y]$.  These notations extend to actions  of a group on
an additive group; see \ref{commdefs}, ff.

Here are some fairly standard 
notations used for particular extensions of
groups: 
$p^k$ means an elementary abelian $p$-group;
$A.B$ means a group extension with normal 
 subgroup $A$ and quotient $B$;  
$p^{a+b+\dots }$ means an iterated group extension, with factors $p^a, p^b,
\dots $ (listed in upward sense); 
$A{:}B, A{\cdot}B$ mean, respectively, a  split extension, nonsplit
extension.


\bigbreak

\section{Introduction}

All lattices in this article are positive definite.  
A sublattice is simply an additive subgroup  of a lattice
(no requirement on the rank).

We prove existence and uniqueness of the Barnes-Wall 
lattices of rank $2^d$ by
induction and establish  properties of them 
and their automorphism
groups, including some new ones.   In particular, the uniqueness
theorem seems to be new.   
With future classifications (and discoveries!) of lattices in mind,
we  promote systematic study of
uniqueness for important lattices.  
In \cite{POE}, we used scaled unimodular lattices and SSD
involutions to give a new uniqueness proof of the Leech lattice and 
revise
the basic theory of the Leech lattice, Conway groups 
and Mathieu groups.  
 There is a new and elementary uniqueness proof for the $E_8$ lattice
in 
\cite{GrE8}.

The Barnes-Wall lattices $\bw d$ are even lattices in
Euclidean space of dimension $2^d$.  They have minimum norm
$2^{\lfloor \frac d 2 \rfloor }$ and remarkable automorphism groups
\cite{BRW2}  isomorphic to $BRW^0(2^d, +) \cong \ratholoex d$, 
$d \ge 4$.    

Various
terms have been applied to these abstract groups and their analogues over
finite fields in general.  We think that {\it BRW group} for the
groups which occur here  would be most appropriate since Bolt, Room
and Wall seem to have been the first to determine 
their structure \cite{BRW2}. 
Compare
the later articles \cite{BE}, \cite{Grex},  \cite{GrDemp},
\cite{GrNW}. See Appendix A2.
 
These lattices (and related ones) were defined in \cite{BW}.  
Independently,  these lattices were rediscovered and 
their groups analyzed by Brou\'e and Enguehard in \cite{BE}.  This
coincidence does not seem well recognized in the literature.  We
first noticed
\cite{BE}, then \cite{BW} only years later.  
  The beautiful and definitive analysis of Brou\'e and
Enguehard 
\cite {BE} was the main inspiration for this article.  

We shall abbreviate Barnes-Wall by BW.

For ranks  $2^d\le 16$, the BW lattices are well-known in
several contexts.  
For $d=1$, we have a square lattice, and, depending on scaling,
$\bw 2$ is the $D_4$ or $F_4$ root  lattice.  We have $\bw 3 \cong
\leh$, though in \cite{BW}, we find $\rtleh$.  As sublattices of many of
the Niemeier lattices, there are scaled copies of $\bw 3$ and $\bw 4$.  
See also
\cite{POE}.  About the BRW groups, there are further details in Appendix
A2.

We prove existence and uniqueness of the BW
lattices of rank $2^d$ by
induction and establish basic properties of them and their automorphism
groups.   
We start 
not with  a frame (a double orthogonal basis) but an orthogonal  sum
of two scaled BW lattices of rank $2^{d-1}$, then show how, by
choosing overlattices, to enlarge this  to  a BW lattice of rank
$2^d$.  Analysis of choices and induction
 give suitable existence and 
uniqueness theorems, structure of the set of minimum norm vectors,  
properties of  automorphism groups, transitivity on certain
sublattices,  etc.  The uniqueness and transitivity theorems are 
new.

 Our program emphasizes elementary algebra and involves 
very little
of special calculations, matrix work and combinatorial arguments.
We heavily exploit commutator density and equivalent properties, 
like
$3/4$-generation and $2/4$-generation,  which are quite useful for
manipulating sublattices and lessening computations.  As far as we
know, these properties are new.

Reflections on the uniqueness theory led us naturally to the {\it
Ypsilanti lattices}, a very large family of BW-like lattices.  The
Ypsilanti lattices are fairly explicit and represent a nontrivial 
share
of all the even unimodular lattices of dimension $2^d$.    Their
existence also clarifies the need for some hypothesis like (e) in
\ref{condx}, as we now explain.  

Let $n>0$ be an integer divisible by 8.  If  $L$ is a rank $n$ even,
unimodular lattice, the theta function of $L$ lies in a vector space
of dimension roughly $n \over 24$ (see \cite{Se}, p.  88).  
For $n=8$ or $16$, the dimension is 1, so the condition constant term 1
determines the theta function.  For $n=24$, the two conditions
constant term 1 and no roots determines the theta function.  In these
cases, one can use arithmetic information about norms to
determine structure.

Now take
$n$ to be $2^d$ for $d$ large and $L$ a BW lattice.  The
condition minimum norm 
$\mu (L) = 2^{\lfloor {d \over 2} \rfloor }$ represents 
$2^{\lfloor {d \over 2} \rfloor }$ 
 linear demands on the theta function.  This number is much
less than $2^d \over 24$.  It is unclear how knowledege of some
higher coefficients can  be used effectively to determine
structure.  The family of Ypsilanti
lattices shows that many isometry types in a given dimension have the
same minimum norm.  To characterize these, or ones like them, we
probably need more than hypotheses about their theta functions.  
We guess that for the Ypsilanti lattices,  
given theta functions may be shared by large
sets of isometry types, and similarly for automorphism  groups.

We acknowledge helpful conversations with  Alex Ryba, 
Leonard Scott, 
Jean-Pierre Serre and Kannan Soundarajan.  

The  author has been supported by
NSA grant USDOD-MDA904-03-1-0098.

\section{Statement of Results}

First, we give some notation, then  state the main results.

\begin{de}\labttr{mu} 
Given a lattice, $L$, 
define $\mu (L):= min\{ (x,x) | x \in L, x \ne 0
\}$.
\end{de} 

\begin{de}\labttr{discgr} Given a lattice $L$, we define the {\it
dual lattice} to be $\dual L :=\{ x \in \QQ \otimes L | (x,L)\le \ZZ
\}$.  Given an integral lattice,
$L$, we define the {\it discriminant group of $L$} to be 
$\dg L := \dual L / L$, a finite abelian group.  A set of {\it
invariants}  of an integral lattice are the orders of the cyclic
summands in a  direct product decomposition of $\dg L$.  (This depends
on choice of decomposition.)  
\end{de}

\begin{de}\labttr{condx} {\bf Condition $X(2^d)$: }
This is defined for integers $d\ge 2$.  Let $s \in \{0,1\}$ be the
remainder of $d+1$ modulo 2. 

We say that the quadruple $(L,L_1,L_2,t)$    
is a {\it an X-quadruple} if it 
satisfies {\it condition $X(2^d)$} (or, more simply, 
{\it condition X}), listed below:

(a)  $L$ is a rank $2^d$ even integral lattice containing 
$L_1\perp L_2$, 
the  orthogonal direct sum of 
sublattices $L_1\cong L_2$ of rank
$2^{d-1}$;

(b) When $d=2$, $L\cong L_{D_4}\cong BW_4$ 
and $L_1\cong
L_2\cong L_{A_1^2}$; when $d\ge 3$, 
$2^{-{s\over 2}}L_1$ and $2^{-{s\over 2}} L_2$ 
are initial entries of quadruples which satisfy
condition  $X(2^{d-1})$.  

(c) $\mu (L) = 2^{\lfloor {d\over 2} \rfloor }$.

(d)   $\dg L
\cong 2^{2^{d-1}}, 1$ as $d$ is even, odd, respectively.  

(e) There is an isometry $t$ of order 2 on $L$ 
which interchanges $L_1$ and $L_2$ and 
satisfies $[L,t]\le L_1\perp
L_2$, i.e., acts trivially on $L/ [L_1 \perp L_2 ]$.  
\end{de}

\begin{de}\labttr{bwt} 
Also, we say that the lattice $L$ is a {\it  lattice of Barnes-Wall type}
or a {\it Barnes-Wall type lattice} if there exist sublattices $L_1, L_2$
of
$L$ and an involution $t \in Aut(L)$ so that $(L,L_1,L_2,t)$   
satisfies condition  $X(2^d)$.  

\end{de}

\begin{thm}\labtt{maintheorem}  Let $d\ge 2$.  A Barnes-Wall type lattice
of rank
$2^d$ exists and is unique up to isometry.
\end{thm}

\begin{coro}\labtt {maincoro}
(i) 
For every integer $d \ge 2$, there is an  integral even lattice
$L$, unique up to scaled isometry, such that

(a) the rank is  $2^d$;  

(b) $Aut(L)$ contains a group $\gd d \cong \ratholoex d $;

(ii) For such a lattice, the group of isometries is isomorphic to
$W_{E_8}$ if
$d=3$ and is just  $\gd d$ for $d \ge 4$.  Also, 
$\dg L \cong 1$ or $2^{2^{d-1}}$, as $d$ is 
odd, even, respectively.   Also, $\mu (L) = 2^{\lfloor {\frac d 2}\rfloor }$.

\end{coro}

We mention that the much-studied lattice $\leh$ is the case $d=3$ of the
above.  The author has recently given an  elementary uniqueness proof for
$\leh$.  See \cite{GrE8},
where previous uniqueness proofs are discussed.  Also, a uniqueness proof
for
$BW_{4}$ was given in
\cite{POE}.

In addition we prove transitivity results for certain types of sublattices
made of scaled Barnes-Wall lattices, including frames.  See
\ref{G3onnorm2},
\ref{transdirty}, \ref{transclean}.   

A final application  of our theory is the construction of {\it the
Ypsilanti lattices} or the {\it Ypsilanti cousins}, built in a similar style.   
(Their definition is a special case of $\frak Y$ 
\ref{fraky}, which is in turn a natural extension of the 
notation ${\frak X}$ \ref{condx};  
the idea came during a pleasant moment in Ypsilanti, Michigan.)  

Let $ j\ge 1,  d=5+3j$.    The Ypsilanti lattices are
indecomposable,  even, unimodular in dimension
$2^d$,  and BW-like in the sense of minimum norm.  For large
dimensions, they become quite numerous.  
 The following
easily stated results give a sample of what we proved.  

\begin{thm}\labtt{manyypsicousins}  
For $c \in [0,\eighth )$ and integer $j>0$ so that 
$\sixteenth (2-2^{1-j}+ 3\cdot 2^{-1-2j}) > c$, 
there is a family 
$Ypsi(2^d,j)$ of rank
$2^d$ indecomposable,  even unimodular lattices, defined for all 
$d>>0$, so that $log_2$ of the number of isometry types in $Ypsi(2^d,j)$ 
has dominant term at least $c \, d \, 2^d$ (in other language, at
least
$(\eighth +o(1))d\, 2^d$.

(i) There is an integer $m$  so that for 
$d>>0$ and 
$L\in
Ypsi(2^d,j)$, 
$\mu (L) = 2^m$.  

(ii) the minimal vectors of $L$ span a proper sublattice 
of finite index in
$L$;

(iii) $Aut(L)$ 
has a normal 2-subgroup $U$ of order divisible by $2^{1+2d}$
\end{thm}

The quotient $Aut(L)/U$ is generally small.  The integer $m$ in (iii)
is roughly $\lfloor {d-j \over 2} \rfloor$.  Like the BW lattices,
the minimum norms go to infinity roughly like the
square root of the dimension.   

\begin{coro}\labtt{maincor} Let $b\in [0,\eighth )$.  
The number $\Upsilon (n)$ 
of isometry types of even unimodular lattices of dimension $n\in 8\ZZ$
which contain a  Ypsilanti lattice as an orthogonal direct summand
satisfies: $log_2(\Upsilon (n))$ is asymptotically at least 
$b\cdot  log_2(mass(n))$ for
$n>>0$, where
$mass(n)$ is the number provided by the Siegel mass formula.  
\end{coro}

\section{Background on Codes}

\begin{de} \labtt{paritycheckmatrix}\rm  
An $(n-k) \times n$ matrix of the form $H= (A|I_{n-k})$,
where
$A$ is an
$(n-k)
\times k$ matrix, is a {\it parity check matrix for the code $C$} if $C$
is defined as the set of row vectors $x \in F^n$ which satisfy $Hx^{tr}=0$ 
\cite{MS}, p.2.  
\end{de}

\begin{de}\labtt {hcsc} \rm 
{\it The Hamming code} ${\cal H}_r$ is defined (up to coordinate
permutations) by the parity check matrix
$H_r$ which is the
$r \times (2^r-1)$ matrix consisting of the $2^r-1$ nonzero column vectors
of height $r$ over $\FF_2$.  
{\it The binary simplex code} ${\cal S}_r$ is the annihilator of the
Hamming code
${\cal H}_r$. 

\end{de} 

\begin{rem}\labtt{hamminginterpretation} \rm The code ${\cal H}_r$ can be
interpreted as the subsets of nonzero vectors in $\FF_2^r$ which sum to
zero.  It has parameters
$[2^r-1,2^r-1-r,3]$ \cite{MS}, p. 23. 
The minimum weight elements of ${\cal H}_r$ are simply the nonzero
elements of a 2-dimensional subspace.  Therefore, a nonzero codeword $A$
in the annihilator meets every such 3-set in 0 or 2 elements. 
Equivalently, the complement $A'$ of $A$ meets every such 3-set in a
1-set or the whole 3-set.  It is clear that $A'$ with the zero vector is
a codimension 1 linear subspace of $\FF_2^r$, whence $A$ is an affine
codimension 1 subspace.  
It follows that every nonzero element
of ${\cal S}_r$  has weight $2^{r-1}$, so  ${\cal S}_r$  has parameters
$[2^r-1,r,2^{r-1}]$.  Note that for $r \ge 2$, ${\cal H}_r \ge {\cal
S}_r$ and that ${\cal S}_r$ contains $\vac$, the all-ones vector, an odd
set.  Also, ${\cal H}_r$ is spanned by the affine planes with 0 removed.  
\end{rem}

\begin{de}\labtt {xhcsc} \rm 
 The {\it extended Hamming code} is obtained by appending
an overall parity check, so has parameters $[2^r,2^r-r-1,4]$.  For $r \ge
2$, it contains the all-ones vector. It is denoted ${\cal H}_r^e$.  Its
annihilator is the {\it extended simplex code} 
${\cal S}_r^e$, which has parameters $[2^r,r+1,2^{r-1}]$.  We have 
 for $r \ge 2$, that ${\cal H}_r^e \ge {\cal S}_r^e$ contains
$\allones{2^r}$. Also, ${\cal H}_r$ is spanned by the affine planes.  

\end{de}

\begin{prop}\labtt{authamm}  If $r\ge 1$ is an integer, the Hamming code
and simplex code of length $2^r$ have automorphism group isomorphic to
$AGL(r,2)$.
\end{prop}
\pf  This is well known.   Since these two codes are mutual annihilators,
they have a common group.   A recent proof was given in an appendix of 
\cite{DGH}.  
\eop

\begin{lem}\labtt{rvl} 
If $S$ is a subset of $\FF_2^d$ of cardinality $2^r>1$ so that for every
affine hyperplane $H$ of $\FF_2^d$, $|H\cap S|=0, 2^r$ or $2^{r-1}$, then
$S$ is an affine subspace.
\end{lem}
\pf 
 This is a result of Rothschild and Van Lint,
\cite{RVL}; it is given in
\cite{MS}, Chapter 13, Section 4, Lemma 6, page 379. 
\eop 

\begin{rem}  \labttr{reedmuller} The Reed-Muller codes are present in our
analysis (the codes ${\cal C}_X$ in \ref{labeling}) but play a small
role.  
\end{rem}  

\begin{de}\labttr{defindeccode}  A code $0 \ne C\le F^X$ is {\it
decomposable} if there is a nontrivial partition $X=Y\cup Z$ of the
index set, so that
 $C=C_Y\oplus C_Z$ is a nontrivial direct sum, where $C_W$ means the set
of vectors in $C$ with support contained in $W \subseteq X$.  If a code
$C\ne 0$ is not decomposable, it is {\it decomposable}.  
\end{de}

\begin{lem} \labtt{indeccode} 
For all $t\ge 3$, there is a length $2^t$ indecomposable doubly even self
orthogonal binary code.  
\end{lem} 
\pf 
For $t=3$, take the extended Hamming code.  Suppose $t \ge 4$ and set
$u=t-3$.  Take a partition of an index set $S$ of size $2^t$ into $2^u$ 
parts $S_i$ of size 8, for $i=1, \dots  , 2^u$.  Let $H_i$ be an extended
Hamming code on
$S_i$.  Take a vector $v_i$ of weight 2 with support $A_i$ in $S_i$ and
define
$v=\sum_i v_i$.   Form the code $C$ spanned by $v_i$ and the codimension
1 subspace of $\sum_i  H_i$ which annihilates $v$.    
Then $wt(v)=2.2^u \in 4\ZZ $, whence $C$ is even.  
\eop

\section{Background on Lattices}

\begin{lem}\labtt {kstlattices}
Let $L$ be a positive definite integral lattice.  Then $L$ has a unique
orthogonal decomposition into indecomposable summands.  More precisely,
let $X(L)$ be the set of nonzero vectors of $L$ which are not expressible
as the orthogonal sum of two nonzero vectors of $L$.  Generate an
equivalence relation on $X(L)$ by relating two elements if their inner
product is nonzero.  An orthogonally indecomposable summand of $L$ is the
sublattice spanned by an equivalence class in $X(L)$.  
In fact, an orthogonal direct summand is a sum of a subset of this set of
sublattices.  
\end{lem}

\pf (See \cite{Kneser} and \cite{MH}, which credits
\cite{Kneser}.)  Let $X_i$, $i=1,\dots , t$ be the equivalence
classes in
$X=X(L)$ and let
$L_i$ be the sublattice spanned by $X_i$.   Positive definiteness implies
that 
$L$ is the sum of the $L_i$.  Also by taking inner products, we deduce $L_i
\cap L_j = 0$ for
$i \ne j$.  So, we have an orthogonal direct sum.  

Let $M$ be an arbitrary orthogonal direct summand.  Let
$N$ be the annihilator of $M$ in $L$. 
   We show for each $i$ that $L_i \le
M$ or $L_i \le N$. 
For $x \in X$, write $x=x_M+x_N$, where $x_M \in M, x_N \in N$. 
Indecomposability implies that one of these components is 0. 

Now, suppose that $M \cap
L_i$ is nonempty.
Then, there exists $x \in X_i$ so that $(x,M) \ne 0$.  The last paragraph
implies that $x\in M$.  We then deduce $X_i \subset M$ and $L_i \le M$. 
If $M \cap L_i=\emptyset$, then $L_i \subset N$.  
\eop

\begin{nota}\labtt {diagonalnotation} \rm 
Given a lattice $L$, the ambient vector space is $\QQ \otimes L$, with
natural extension of the symmetric bilinear form on $L$. 

Take isometries $\psi_i : \QQ \otimes L \rightarrow V_i$ of rational
vector spaces.  From these, we get isometries 
$\psi_{ij}=\psi_i^{-1} \psi_j$ from $V_i$ to $V_j$.  Priming on an index
means replacement of the corresponding map by its negative.

 For a subset $A \subseteq  \QQ \otimes L$, define the
following subsets of
$V := V_1 \perp V_2$: 

$A_{ij}:=\{ x^{\psi_i} + x^{\psi_j} | x \in A \}$,

$A_{ij'}:=\{ x^{\psi_i} - x^{\psi_j} | x \in A \}$,

$A_{i'j}:=\{ -x^{\psi_i} + x^{\psi_j} | x \in A \}$,

$A_{i'j'}:=\{ -x^{\psi_i} - x^{\psi_j} | x \in A \}$. 
\end{nota}

\begin{nota}\labtt{lbc}\rm
Given a basis 
${\cal B}$ of Euclidean space and binary
code 
$C \le \FF_2^{\cal B}$, we define 
$L_{{\cal B},C} := 
\{ \sum_{b \in {\cal B}} \half a_b b | 
a_b \in \ZZ,  (a_b + 2\ZZ)_{b \in {\cal B}} \in C \}$.   (This lattice is
sometimes integral.) 
Note that 
$L_{{\cal B},C} = \{ \sum_{b \in {\cal B}} \half a_b b | a_b \in \ZZ, 
\sum_{b \in {\cal B}} (a_b+2\ZZ) c_b = 0 +2\ZZ \hbox{ for all }
(c_b)_{b \in {\cal B}} \in C^\perp  \}$.   
\end{nota}

\begin{nota}\labtt {twistede8}\rm 
Let $\a_i$, $i \in \Omega = \FF_2^3$ be vectors in $\RR^\Omega$ which
satisfy $(\a_i, \a_j )=2\kron ij$.  
Let $\hame 8$ be the extended Hamming code \ref{xhcsc}.

Define the {\it $A_1^8$-description of $\leh$} or the {\it 2-twisted
version of $\leh$} to be the $\ZZ$-span of all $\a_i$ and all $\half
\a_c$, for $c
\in
\hame 8$.   In the Notation \ref{lbc}, this is $L_{\{\a_i|i =1,\dots ,
8\},
\hame 8}$.
\end{nota}

\begin{nota} \labtt{sumsandsigns}  \rm
Suppose that $\Omega$ is an index set and $\{v_i|i \in \Omega \}$ is a
basis of a vector space.  
For a subset $A$ of $\Omega$, define $v_A:= \sum _{i \in A} v_i$.  The
linear transformation   
$\varepsilon _A$ sends $v_i$ to $-v_i$ if $i  \in A$ and to $v_i$ if 
$i \not \in A$.  The group of such maps is ${\cal E}_{P(\Omega )}$.  If $C$
is a subset of $P(\Omega )$, 
${\cal E}_C$ denotes the set of maps $\varepsilon_A$  for $A \in C$.  This
is a subgroup if $C$ is a subspace of the vector space $P(\Omega )$.  
\end{nota}

\begin{prop}\labtt {lbhame} 
For an integer $d \ge 2$, define $m:=\lfloor {d \over 2}\rfloor$.  Let
$\Omega$ be an index set, identified with $\FF_2^d$.  
Take a basis
${\cal B}:=\{ v_i | i\in \Omega \}$ where $(v_i,v_j)=2^m \kron ij$ of
$\RR^\Omega$.  Form $L_{{\cal B, H}^e}$, as in \ref{lbc}; it is
integral for $d\ge 3$.   
Then, if $d \ge 4$, $Aut(L_{{\cal B}, \hame d})$
is in the monomial group on $\cal B$ and in fact 
$Aut(L_{{\cal B}, \hame d}) = {\cal E}_{\Omega}F$, where $F$ is a natural
$AGL(d,2)$ subgroup of the group of permutation matrices.  If $d=3$, 
$Aut(L_{{\cal B}, \hame d}) \cong \weh$.  
\end{prop} 

\pf  For $d=3$, we have the lattice \ref{twistede8} and for $d=2$, we have
the $F_4$ lattice, spanned by vectors of shape $(\pm 1,0,0,0), 
(\pm \half, \pm \half, \pm \half, \pm \half)$.  These automorphism groups
are well known to be $\weh$ and $W_{F_4}$, respectively.

 For any $d\ge 2$, the set  of minimal vectors of 
$L:=L_{{\cal B}, \hame d}$ is just $\pm v_i$,
$i\in \Omega$ and $\half v_S \varepsilon_T$, for $S \in \hame d$ an
affine plane and $T$ a subset of $S$.  
These span $L$ since affine
planes span
$\hame d$.   All these minimal vectors have norm $2^m$, $m \ge 1$.  

We now  assume
$d
\ge 4$.  
The
set of these which are in
$2^{m-1}L^*$ is exactly $\pm v_i$, $i \in
\Omega$, for if $A$ is an affine plane there exists an affine plane $B$
so that $A\cap B$ is a 1-set (because $d \ge 4$), whence $(\half  v_A
\varepsilon_S,
\half v_B
\varepsilon_T )=\pm 2^{m-2}$.  
It follows that $Aut(L)$ is contained in 
the monomial group based on $\cal B$.  Clearly it contains  ${\cal
E}_{\Omega}F$, described above, and maps to the stabilizer of the code
$L/Q$ in $\half Q/Q\cong \FF_2 ^ \Omega$, where $Q$ is the square lattice
with basis
$\cal B$.  Since $Aut(\hame d )$ is a natural $AGL(d,2)$ subgroup of the
symmetric group (\ref{authamm}), we are done.  
\eop

\begin{de}\labttr{euclift}  Let $c=(c_i)\in \FF_2^n$.  
The {\it Euclidean 
lift} of $c$ is the vector in $\{0, 1\}^n \subset \ZZ^n$ 
which reduces modulo 2 to $c$.
When $p$ is an odd prime and $c=(c_i)\in \FF_p^n$, we have a similar
definition of lift, using   the subset $\{-{p-1\over 2}, -{p-3\over
2} \dots ,-1, 0, 1, \dots {p-3\over 2}, {p-1\over
2} \}^n
\subset
\ZZ^n$.  
\end{de} 

\begin{lem} \labtt{projweight} 
Let $L$ be a lattice with sublattice of finite index $M$ which is a 
coelementary abelian $p$-group for some prime $p$.  
Let $F:=\FF_p$.  
Suppose that 
 $\cal C$ is an error correcting code in $F^n$ with minimum weight
$w$.  
Suppose that $J$ is the lattice between 
$M^n$ and $L^n$
corresponding to $\cal C$, i.e. spanned by all 
$(c_1x, c_2x, \dots , c_n x)$ 
for $x \in L$ and $(c_i)$ is the Euclidean
lift of a codeword in $\cal C$.  

(i) 
If $(y_1, \dots , y_n) \in J
\setminus M^n$, the weight of $(y_1, \dots , y_n)(mod \ M^n)$ 
is at least
$w$.  

(ii) Suppose that $M$ is indecomposable and $\cal C$ is
indecomposable.  Then
$Aut(J)\cap Aut(M\perp
\dots \perp M)\cap Aut(L\perp \dots \perp L)$ factorizes as the
product of subgroups $A_1A_2$, where
$A_1$ is the subgroup which fixes each direct summand isometric to
$M$ and acts diagonally on $(L/M)^n$, and where $A_2$ is the subgroup
of the natural group of block permutation matrices of degree $n$ which
fixes the code $\cal C$.   
\end{lem}  
\pf  
(i) 
We take a basis $v(1), \dots , v(d)$ of $\cal C$.  For a codeword
$v=(v_i)$ and  $x \in L/M$, let $vx$ be the vector in
$(L/M)^n$ whose $i^{th}$ component is the Euclidean lift of $v_i$ times
$x + M$.  
So, $J/M^n = \sum_{i=1,\dots ,d; x \in L} v(i)x + M^n= \bigoplus
_{i=1,\dots ,d}v(i)(L/M)$.

Suppose that  $y=(y_1, \dots , y_n)(mod \ M^n)$ represents the minimal
weight $k$ in $J/M^n$.  Write it  as a linear combination 
$y=\sum_{i=1}^d v(i)x(i)$, where $x(1), \dots, x(d)$ is a sequence of
elements of $L/M$ and the product $v(i)x(i)$ is as in the previous
paragraph.

Take any linear functional $f:L/K\rightarrow F$ and extend it
componentwise to
$g:(L/M)^n
\rightarrow F^n$.  Then $g(J/M)={\cal C}$.   Given a nonzero $x(i)$,
take an
$f$ so that
$f(x(i))=1$. Then $g(y) \in {\cal C}$ has nonzero coefficient 1 at $v(i)$,
whence
$0\ne g(y)=(f(y_1), \dots , f(y_n)) \in {\cal C}$ has at least $w$
nonzero entries, whence so does $y=(y_1 , \dots , y_n )$. 
Therefore,
$k
\ge w$.   

(ii)  
First, suppose that 
$h\in Aut(J)\cap Aut(M\perp \dots \perp M) \cap 
 Aut(L\perp \dots \perp L)$. 
We claim that $h$ determines a unique element of the code group, 
up to scalars.  
For any $v \in  {\cal C}, x \in L \setminus M$, $h$ takes $vx$ to an
element of  
$J/M$.  
This means that for  
linear functionals $f, g$ as in (i) where $f(x) = 1$, 
we have that $g(h(vx))$ is a codeword. 
Since $h$ takes $vx$ to another element of similar form $v'x'$, it 
follows that there is a block permutation matrix $b$ so that, for all
codewords $x$,  the action of $hb$ stabilizes each direct summand
isometric to 
$M$ and takes
$vx$ to an element of the form
$vx'$, for all $x \in L/M$.   Since the code is indecomposable, we
use the property that for any $v, w \in {\cal C}$, if $hb$ takes $vx$
to $vx'$, then $hb$ takes $wx$ to $wx'$.  In other words, the actions
of  $hb$ on the summands of $(L/M)^n$ are identified.  
\eop

\begin{lem}\labtt{indecoversmv} Suppose that $J$ is a lattice and that
$SMV(J)$, the sublattice spanned by the minimal vectors, has finite
index.  Suppose that
$SMV(J)=J_1\perp
\dots 
\perp J_n$, where the $J_i$ are indecomposable lattices.  

Let $K$ satisfy $SMV(J) \le K \le J$ and $K$ is homogeneous  with
respect to the rational vector spaces spanned by the summands, i.e., 
$K=\sum_{i=1}^n K_i$ where $K_i:= K \cap (\QQ \otimes J_i)$.

Suppose that
$J/K$ corresponds to an indecomposable code \ref{defindeccode}.  Then $J$
is orthogonally indecomposable.  
\end{lem}  
\pf (See \ref{kstlattices}.)  
Let $x$ be an indecomposable vector of $J$ which is not in $K$ and let $S$
be the indecomposable summand of $J$ which contains it.  Let $A$ be the
support of $x+K$ in $J/K$, i.e., those indices where $x+K$
projects nontrivially to $\QQ \otimes K_i/K_i$.  For 
$i
\in A$, there exists a minimal vector $y \in K_i$ so that $(x,y)\ne 0$.  
Therefore,
$J_i \le S$, for all $i \in A$.  The indecomposability assumption on the
code implies that all $J_i, i=1, \dots ,n$ are in $S$ and since
$SMV(J)$ has finite index in $J$ and 
$S$ is a summand of
$J$,
$S=J$.  
\eop

\begin{de}\labtt {rmodlatt} \rm  Let $r>0$ be an integer.  
An integral
lattice
$L$ is {\it $r$-modular} if $L \cong \sqrt r \dual L$ 
\end{de}  

\begin{de}\labtt{ssd}\rm   The SSD concepts were established in
\cite{POE}.  
Call a lattice
$M$ {\it semiselfdual} (SSD) if
$2\dual M \le M\le \dual M$.  If the sublattice $M$ of the integral
lattice $L$ is semiselfdual, we define the orthogonal transformation 
$t_M$ on $\QQ\otimes
L$  by $-1$ on $M$ and $1$ on $M^\perp$.  Then $t_M$
leaves $L$  invariant and so gives an isometry of $L$ of order 1 or 2;
it has  order 2 on
$L$ if $M
\ne 0$. 

A more general notion is that of {\it relatively SSD (RSSD)}: this is the
condtion that the sublattice $M$ of the integral lattice $L$ satisfies
the weaker condition 
$2L \le M +  M^\perp $.    In this case, the orthogonal involution defined
as above preserves $L$.  
\end {de}

\section{Actions of 2-groups and endomorphisms on lattices}

We gather an assortment of results on this topic.

\begin{de}\labtt{fourvolution} 
\rm A {\it fourvolution} is a linear
transformation whose square is $-1$.  
A {\it fourvolution on a lattice} is a lattice isometry whose square is
$-1$.  In case we have a lower group as in \ref{lowerstuff}, we use
the terms {\it lower and upper fourvolution}.  
We may call an element in a group a fourvolution with respect to 
a representation, and even with respect to more than one representation, for 
example by restriction of one action  to a submodule.  
\end{de}

\begin{lem}\labtt {aboutfourvolutions}  If $f$ is a fourvolution of the
lattice
$L$, then the adjoint of $1\pm f$ is $1 \mp f$, 
$(1\pm f)^2=\pm 2f$, $1\pm f$ is an isometry scaled by $\sqrt 2$ and we
have
$L \ge L(1+f) \ge 2L$ and $|L:L(1+f)|=|L(1+f):2L|=|L/2L|^{\half}$.  In
particular,
$rank(L)$ is even.   
\end{lem} 
\pf 
For $x, y \in L$, we have $(x(1\pm f),y(1\pm f))=(x,y) \pm (x,yf) \pm
(xf,y) \pm (xf,yf)= 2(x,y)\pm (x,yf)\pm (xf^2,yf)=2(x,y)$.  For
adjointness, just compute $(x(1 \pm f),y)=(x,y)\pm (xf,y)= (x,y)\pm
(xf^2,yf)=(x,y)\mp (x,yf)$.    The other  statements are easy to prove.  
\eop

\begin{nota}\labtt {twist}  \rm 
Let $L$ be a lattice
and $f$  a
fourvolution in $Aut(L)$.  Define
$S[k] := S(1-f)^k$, for
$S
\subseteq
\QQ
\otimes L$ and $k \in \ZZ$.  Note that this makes sense since the linear
map
$1-f$ is invertible, with inverse 
$\half (1+f)$.  Call a transformation  of the form $1-f$ a {\it sultry
tranformation}
and
call 
$S[k]$ the {\it $k^{th}$ sultry $(1-f)$-twist of $S$} or
the {\it $k^{th}$ sultry twist of $S$}.  
(The terminology is explained in \ref{sultryinterpretation}.)

  We have $S[0]=S$.  Note that for all $k$, $(Sp)[k]=(S[k])p$, where $p$
is any polynomial expression in $f$.  Also, $S[k][\ell ] = S[k+\ell ]$.

\end{nota}

\begin{lem}\labtt{dualityleveltwist} If $S$ is an $f$-invariant 
lattice in
$\QQ
\otimes L$, then for $k \le \ell$, $|S[k]:S[\ell ]|=2^{\half
rank(S)(\ell - k )}$.  
\end{lem}\pf
This follows since $(1-f)^2=-2f$ and because for all integers $p, q$ and
all integers $r\ge 0$, $S(1-f)^p/S(1-f)^{p+r} \cong S(1-f)^q/S(1-f)^{q+r}$.
\eop

\begin{lem}\labtt {dottwists} 
Let $S$, $T$ be subsets of $\CC \otimes L$.  Then 

(i) $(S[1],T)=-(S,Tf[1])$.

Now assume that $S$ and $T$ are $f$-invariant. Then $S=Sf=-S$, $T=Tf=-T$
and the following hold.

(ii) For  all integers $k, \ell$, we have 
$(S[k],T[\ell ])=2(S[k-1],T[\ell -1])$ and $(S[k],T[\ell
])=2(S[k-2],T[\ell])=2(S[k],T[\ell -2])$.

(iii) $S^*[k]=S[k]^*f^{-k}$.

(iv) $(S[k],T[\ell])=(S[k'],T[\ell'])$, for all integers $k, k', \ell,
\ell'$ such that $k+\ell = k' + \ell'$; and 

(v) Assume that the integer 
$\ell$ satisfies $S^*=S[\ell ]$.  Then
$S^*[k]=S[k+\ell ]$.   

\end{lem} 

\pf (i) and (v) are clear. 

(ii) follows since $1-f$ is an isometry scaled by $\sqrt 2$. 

(iii) We have $x \in S[k]^*$ if and only if $(x,S[k]) \in \ZZ$ if and only
if $x(1+f)^{k}=(-1)^kxf^k(1-f)^k \in S^*$ if and only if $x \in
(-1)^kS^*[-k]f^{-k} = S^*[-k]f^{-k}$.  

(iv)  is trivial 
for $k=0$  and for $k \ge 1$ it follows from (i) and easy
induction.  If $k$ is negative, use (ii) and the case $k \ge 0$.

\begin{ex}\labtt{rank4twist} \rm If  $L\cong L_{D_4}$ then $L[-1]
\cong L_{F_4}$, where we take the latter to be the span of a standard 
version of the $F_4$ root system: $(\pm 10^4), (\pm 1^20^2), (\pm \half
^4)$. 
\end{ex}

\begin{de}\labtt {dualitylevel} \rm 
Let $L$ be a lattice with fourvolution $f$.  
Suppose that there is an integer $r$ such that 
$L^*=L[-r]$ (see \ref {twist}). 
We call $r$ the  {\it duality level} of $L$.  Such a modular lattice
(see \ref{rmodlatt} and \ref{dualityleveltwist}) is called an {\it
$r$-sultrified dual} and is $2^r$-modular (see \ref{rmodlatt}).  
\end{de}

\begin{de}\labtt {eigenlattice} \rm 
Given a group $E$ acting on a lattice $L$ and character $\l \in
Hom(E, \{ \pm 1 \} )$, define the {\it eigenlattice } $L^\l$ to be $\{ a
\in L | ay=\l (y) a, \hbox{ for all } y \in E\}$.  Define the {\it total
eigenlattice} to be $Tel(E,L):=\sum_{\l \in Hom(E, \{ \pm 1 \} )}
L^\l$.   The notation extends naturally a set of
automorphisms.  When $t$ has order 1 or 2, define $L^+, L^-$ to be the
lattice of fixed,  negated points, respectively.  To denote dependence
on
$t$, we write $L(\pm,t)$ or $L^\pm(t)$ for $L^\pm$.  
\end{de}

\begin{rem} \labtt{noncyclicelab} \rm  
In case $E$ is 2-elementary
abelian,
$L/Tel(E,L)$ 
is finite, and is in fact a 2-group, but in
general is not elementary abelian.   For an example, let $E$ be a
fourgroup and
$L=\ZZ [E]$, the regular representation.  Then,  $L/Tel(E,L) \cong 2
\times 2 \times 4$. 
\end{rem}

\begin{lem}\labtt {defectineq}
Suppose that the involution $t$ acts on the additive group $A$.
Let $A^\varepsilon := \{ a \in A | a^t = \varepsilon a \}$.  
Suppose furthermore that 
the minimal number of generators of $A$ as an abelian
group is $r < \infty$. 
Define integers $k, \ell$ by $2^k := | A{:} A^- + A^+|$ and $2^\ell
:=|A{:}B|$, where $B:= \{ x \in A | x(1-t) \in 2A \}$.  

Then: (i) $2A \le A^- + A^+ \le B$, whence $\ell \le k$; 

(ii) 
Then $\ell  \le r/2$.

(iii) In the notation of (ii), $A^- \ge A(1-t) \ge 2A^-$ and 
$|A(1-t)/2A^-|=2^k$, whence $rank(A^-) \ge k$.  

(iv) In the notation of (ii), $A^+ \ge A(1+t) \ge 2A^+$ and 
$|A(1+t)/2A^+|=2^k$, whence $rank(A^+) \ge k$.  

(v) If $A$ is free abelian, $A^\varepsilon$ is a direct summand of $A$
and $A^+ + A^- = A^+ \oplus A^-$. 

(vi) If $A$ is free abelian, $k=\ell$ (whence $k \le r/2$).  

(vii) Suppose that multiplication by 2 is a monomorphism of $A$ (e.g.,
$A$ is free abelian).  If
$k=0$ (i.e., if
$t$ is trivial on
$A/2A$),
$A=A^+ + A^-$.  
\end{lem}

\pf
(i) The proof follows from the equation $2a=(a+a^t)+(a-a^t)$.  

(ii)  Let $B:= \{ x \in A | x(1-t) \in 2A \}$.  Then 
the map $(1-t)$ induces an injection of $A/B \cong 2^\ell$ into $B/2A$,
so  in particular  $\ell \le k$.  
If $x_1, x_2, \dots , x_\ell \in A$  form a basis modulo $B$, then 
$x_1, x_1^t,  x_2, x_2^t,  \dots , x_\ell, x_\ell^t$ 
are independent
modulo
$2A$.  Therefore, $2\ell \le r$.

(iii) For the first statement, notice that the kernel of the map $\phi : A 
\rightarrow A^-/2A^-$, $x \mapsto x(1-t)$ is $A^+ \oplus A^-$ and then use
$Im(\phi)  \cong A/Ker(\phi)$, which has rank $k$.

(iv) This follows from (iii) by replacing  $t$
with  $-t$.

(v) Clearly, $A/A^\varepsilon$ is torsionfree.  The second statement
follows from $A^+ \cap A^-=0$.  

(vi)  This follows from the general classification of free
abelian groups which are modules for cyclic groups of prime order, e.g.,
\ref{indec2}; (74.3) in  
\cite {CR}.  (The result for a cyclic group of order 2 is easy to
prove directly.)  It states that such a module $A$ has the form $F_1
\oplus
\dots \oplus F_p \oplus E_1 \oplus  \dots  \oplus E_q$, where each
$F_i$ is a copy of the regular representation $\ZZ \la t \ra$ and
where each
$E_j$ is infinite cyclic.  By reducing such a decomposition modulo 2,
one deduces that
$k=\ell$. 

(vii)  This is easy to prove directly (of course it is a consequence of
the nontrivial result mentioned in (vi)).  Suppose that the involution
$t$ is trivial on $A/2A$.  Then $t=1+2S$ for some $S \in End(A)$. From
$1=t^2=1+4(S+S^2)$, we deduce that $S+S^2=0$.  For $a \in A$,
$a=a(1+S)-aS$.  One checks that $aS \in A^-$ and $a(1+S) \in A^+$.  
\eop

\begin{de}\labtt {defect} \rm 
The {\it defect} of the involution $t$ acting on the free abelian
group $A$ is  the integer
$k=\ell$, as in 
\ref {defectineq}.  It is the number of nontrivial Jordan blocks for the
action of $t$ on
$A/2A$. 
\end{de}


\def\Lp{L^+}
\def\Lm{L^-}
\def\Lpm{L^\pm}
\def\Lmp{L^\mp}

\begin{lem}\labtt {involonunimodular} 
Let $L$ be a unimodular lattice  and $t$ an involution acting
on $L$.  Then the eigenlattices $L^\varepsilon := \{ x \in L |
 x^t=\varepsilon x \}$ satisfy $\dg {L^+} \cong \dg {L^-} \cong 2^k$,
where
$k$ is the defect of $t$ in the sense of Definition \ref {defect}.  
\end{lem}

\pf
Since each $L^\varepsilon$ is a direct summand of $L$, which is
unimodular, the orthogonal projection takes $L$  onto
$[L^\varepsilon ]^*$.  
The kernel of the map from $L$ to 
$\dual {[L^\varepsilon]} / L^\varepsilon$ is $L^\varepsilon +
L^{-\varepsilon}$, so  from Lemma \ref {defectineq}, we deduce that the
image is elementary abelian, of order $2^k$.   
\eop

\begin{rem}\labtt{involtossd}\rm  
The notion of SSD involution is essentially the same as that of an
involution on a lattice.   Let $L$ be a lattice.  
An involution $t\in Aut(L)$ creates a pair
of eigenlattices, $L^\pm$.   
Since 
 $L^+\perp L^-$ is 2-coelementary abelian in $L$ and $t$ acts
trivially on $\half [L^+\perp L^-]/[L^+\perp L^-]$,
$t$ is an SSD involution  which preserves $L$ (see  \ref{defectineq},
\ref{ssd} and
\cite{POE}). 
\end{rem}

\begin{de}\labtt {scalarsg} \rm  For a group $G$ acting on the
$RG$-module $M$, where $R$ is a commutative ring, the {\it scalar
subgroup} is $$Scalar(G,M) := \{ g \in G | g \hbox{ acts on $M$ as
multiplication by an element of $R^\times$} \}.$$

When $M$ is a free abelian
group, this is just the subgroup of group elements which act as $\pm 1$.  
\end{de}

\begin{de} \labttr{plainframe}  A {\it frame} or {\it plain frame} in a
rank $n$ lattice  is a set of $2n$ vectors of common norm, two of which are
linearly dependent or orthogonal.  
\end{de} 

Later in \ref{sf}, we work with a special case of this.

\subsection{Commutator density, $3/4$-generation and $2/4$-generation} 

The concepts \ref{commutatordense}, \ref{3/4generation}  and results in
this section  seem to be new. 
Commutator density is an unusual property which is very useful for
controlling commutators of an extraspecial 2-group acting on a lattice.

Note that we will be mixing additive and multiplicative commutator
notation.

\begin{de}\labtt {commdefs}\rm 
We recall a few defintions involving groups and modules.  
Let $Q$ be a (multiplicative) group and $S$ a subset of $Q$.  
For $s, t \in Q$, as usual $[s,t]=s^{-1}t^{-1}st$. 
For a module $M$,
$[M,S]$ is the {\it commutator submodule}, 
meaning as usual the additive group spanned by all commutators 
$[x,s]=x(s-1)$,
$x \in M, s\in S$.  
Higher commutators are interpreted by extending these definitions, for
example
$[x,s,t]=x(s-1)(t-1)$,  $[s,t,x]=-[x,[s,t]]$ and $[t,x,s]=-[x,t,s]$.  
\end{de}

\begin{de}\labtt {commutatordense}\rm 
Let $Q$ be a group and $S$ a subset of $Q$.  

{\it $S$-CD}: 
A module $M$ for $Q$ is
called {\it $S$-commutator dense} if 
$[M,Q]=[M,S]$.  
(When $S=\{f\}$, a single element, every element of $[M,Q]=M(f-1)$ is a
commutator.) 

{\it $S$-$k$CD}:
As is common, for the natural number $k$, we use the notation $[M,Q;k]$
for $[M,Q,Q,\dots ,Q]$ ($k$ times).  We say that $M$ is {\it degree $k$ 
$S$-commutator dense} or if
$[M,Q;k]=[M,S;k]$.

{\it $S$-HCD}:
If $M$ has such properties for all $k \ge 1$, we say that $M$ is
{\it  $S$-higher commutator dense}.

{\it $S$-TCD}:
A module $M$ is {\it  $S$-commutator dense on submodules} if all its
submodules are $S$-commutator dense.  In this spirit, we define 
{\it degree $k$ CD} and {\it $S$-HCD} {\it on submodules}.  

When the set $S$ is understood, we may drop $S$ from the preceeding
notations.  Note that commutator density is inherited by quotient modules
but may not be for submodules.  

\end{de}

\begin{lem}\labtt {3sgapp} 
Suppose that the group $Q$ acts on the $\ZZ Q$-module $L$ and that $L$ is
$f$-commutator dense for a fourvolution $f \in Q$ such that 
$[Q,f]$ is scalar on $L$.   Then $[L,Q;k]=[L,f;k]$, for all $k \ge 1$,
i.e.,
$L$ is
$f$-HCD.    In fact,  $[L,Q;k]=[L,f;k]=2^{k \over 2}L$ if $k$ is even, 
and 
 $[L,Q;k]=[L,f;k]=2^{k-1 \over 2}[L,f]$ if $k$ is odd.   
\end{lem}
\pf  
We have $[L,Q]=[L,f]$, whence $[L,Q,f]=[L,f,f]=2Lf=2L$.  Also,
$[Q,f,L]\le [Scalar(Q,L),L]\le 2L$.  The Three Subgroups  Lemma
\cite{Gor, Hup} implies that
$[f,L,Q]\le 2L$, or $[L,Q,Q]\le 2L$, which is $[L,f,f]$.  
The statements $[L,Q;k]=[L,f;k]$, for all $k \ge 1$ are proven 
by induction.  
\eop

\begin{lem}\labtt{dih8onlatt}  Suppose that the lattice $L$ contains the
orthogonal direct sum of sublattices $L_1\perp L_2$, that $L_1$ and
$L_2$ have rank $n=2m\in 2\ZZ$ and
$L/L_1\perp L_2$ is elementary abelian of order $2^m$.  
Suppose that involutions $t, u$ act on $L$ so that $L_1=L^-(t)$,
$L_2=L^+(t)$ (see \ref{ssd}) and $u$ interchanges $L_1$ and $L_2$. 
Then: 

(i) $u$   acts trivially on
$L/L_1\perp L_2$ if and only if $det(L^\pm (u))=det (L_1)=det(L_2)$.  

(ii) If the conditions of (i) hold, then 
$L$ is the sum of any three of the four sublattices 
$L^\pm (t), L^\pm (u)$.  
\end{lem}
\pf 
Clearly $u$ acts on $\half L_1 \perp \half L_2$ and on 
 $\half L_1 \perp \half L_2/ L_1 \perp L_2$, it acts with $n$ Jordan
blocks of size 2.  
Also, $det((L_1 \perp L_2)^\pm  (u))= det(L_1)2^{n}$ and $det(\half (L_1
\perp L_2)^\pm  (u))= det(L_1)2^{-n}$.  

(i) The equivalence of the two
conditions follows from  comparision of the determinants of the lattices 
$L^\pm (u) \ge (L_1\perp L_2 )^\pm (u)$.  Let $2^r$ be the index 
$|L^\pm (u) {:} (L_1 \perp L_2)^\pm  (u)|$.  Then $det(L^\pm (u) ) =
2^{n-2r}det(L_1)$.  The second condition in (i) implies that $n=2r$,
whence $L^\pm (u)$ covers $L/[L_1 \perp L_2]$.  Conversely, take $x
\in L$.  It is fixed by $u$ modulo $L_1 \perp L_2$, which is a free
module.  Therefore there is $y \in L_1$ with $x(u-1)=y(u-1)$.  The
coset $x+[L_1\perp L_2]$ therefore contains the fixed point $x-y$.  

(ii) 
The hypotheses imply that $L_1+L_2+L^+(u)= L$, and a similar statement
applies to $-u$.  Finally, we may interchange the roles of $t$ and $u$ to
deduce the remaining statements.  
\eop 

\begin{de}\labttr{3/4generation}  Let the dihedral group $D$ of order 8
be generated by involutions $t, u$.  An action of $D$ on the abelian
group $L$ has the {\it 3/4 generation property} if the central
involution of $D$ acts as $-1$ on $L$ and $L$ is the sum of any three of
$L^\pm (t), L^\pm (u)$.  An action has the {\it $2/4$-generation
property} if $L$ is the sum of the fixed points of a pair of generating
involutions.  
\end{de} 

\begin{prop} \labtt{equivcd}   Suppose that the dihedral group $D$ acts
on the lattice
$L$ with the central involution acting as $-1$.  For this action,
equivalent are 
the properties of  3/4-generation, 2/4-generation
and  commutator density for a fourvolution in $D$.  
\end{prop}
\pf  Let $f$ be an element of order 4  in $D$ and $t, u$ a pair of
generating involutions.  Set $L_1:=L^+(t), L_2:=L^-(t)$.  Note that
$rank(L)$ is even.  

Assume the 3/4 generation property, 
and assume the notations of
\ref{3/4generation}.  
Then $L$ has even rank $2n$ and $[L,D] \le
(L_1+L_2)\cap (L^+ (u)
\perp L^-(u))$, which has index $2^n$ in $L$.  
Since $(f-1)^2=-2f$, 
this intersection equals $L(f-1)$, whence density. 

Assume density.  Consider the action of $D$ on $\half (L_1 \perp L_2)/L_1
\perp L_2$.  The action of $t$ is trivial and $f$ acts as an involution
with $n$ Jordan blocks (as $f^2=-1=[t,u]$).  

We have $L > L(f-1)\ge L(t-1)+L(u-1)+L(t+1)+L(u+1)$.  Since 
$\pm t, \pm u$ is a normal subset of generators of $D$, the right side is
$[L,Q]$ which by density equals $L(f-1)$.  Note that $L(t+1)+L(t-1)\ge
2L$ and  $L(u-1)+2L=L(u+1)+2L$.  It follows that 
$L(f-1)=L(t-1)+L(u-1)+L(t+1) \le L^-(t)+L^-(u)+L^+(u)$.

Similar arguments apply if we replace $t, -t, u$ by any 3-subset of $\{t,
-t, u, -u \}$.  This completes the proof  that density implies
3/4-generation. 

Obviously, 2/4-generation implies 3/4-generation.  Assume
3/4-generation and let $t, u$ be any generating pair of involutions.  Set
$M:=L^+(t)+L^+(u)$, a sublattice of $L$.  Since the central involution of
$D$ acts as $-1$, the summands meet trivially, whence $M$ has rank
$2n$.  Since $L^+(t)$ is RSSD in $L$ (see \ref{ssd}), it is RSSD in
$M$, i.e.
$M$ is $t$-invariant.  It follows that $M$ contains $L^-(u)=L^+(u^t)$,
whence $M=L$ by the 3/4-generation property.  
\eop  

We can actually drop reference to the quadratic form in the previous
result.

\begin{prop} \labtt{equivcd2}   Suppose that the dihedral group $D$ acts
on the free abelian group  
$L$ with the central involution acting as $-1$.  For this action,
equivalent are
the properties of  3/4-generation, 2/4-generation
and  commutator density for a fourvolution in $D$.
\end{prop}
\pf  
This follows from  \ref{equivcd} once we define a $D$-invariant 
positive definite integer
valued quadratic form.  One uses the familiar trick of taking any 
integer valued positive definite quadratic form on $L$, then summing its
transforms under $D$.  
\eop

\section{Sultry twists and the NextBW procedure}


\def\lt#1#2{L[#1]_{#2}}
\def\ltt#1#2#3{L[#1]_{#2#3}}

\def\li#1{L_i[#1]} 
\def\lj#1{L_j[#1]} 
\def\lij#1{L_{ij}[#1]} 
\def\lijp#1{L_{ij'}[#1]} 

\def\lone#1{L_1[#1]} 
\def\ltwo#1{L_2[#1]} 

\def\lonetwo#1{L_{12}[#1]} 
\def\lonetwop#1{L_{12'}[#1]} 

\def\dem{M_1[1-r]+M_2[1-r]+M_{12}[-r]} 
\def\demt{M_1[1-r]+M_{12}[-r]+M_{12'}[-r]}

\def\mt#1#2{M[#1]_{#2}}
\def\mtt#1#2#3{M[#1]_{#2#3}}

\def\mi#1{M_i[#1]} 
\def\mj#1{M_j[#1]} 
\def\mij#1{M_{ij}[#1]} 
\def\mijp#1{M_{ij'}[#1]} 

\def\mone#1{M_1[#1]} 
\def\mtwo#1{M_2[#1]} 

\def\monetwo#1{M_{12}[#1]} 
\def\monetwop#1{M_{12'}[#1]} 

\def\del{M_1[1-r]+M_2[1-r]+M_{12}[-r]} 
\def\delt{M_1[1-r]+M_{12}[-r]+M_{12'}[-r]} 

\def\dedl{M_1[1-r]+M_2[1-r]+M_{12}[-1]} 
\def\dedlt{M_1[1-r]+M_{12'}[-1]+M_{12}[-1]}


We discuss some procedures for proving the main theorem.  We continue to
let BW abbreviate ``Barnes-Wall''.  

We first show how
to start from a BW-type  lattice of rank $2^{d-1}$ and 
create one of rank $2^d$. 
Later, in \ref{bwtclassification} 
we show how a BW-type lattice of rank  
$2^d$ is uniquely determined by an ancestor of rank $2^{d-1}$.    
Eventually, we use an induction argument which will show that a 
BW-type lattice is unique, so is the same (up to rescaling) 
as the  
lattices constructed in \cite{BW, BE}.

An important technique here is to use the commutator density enjoyed
by these lattices.  The twisting by sultry transformations helps
control the analysis.

\begin{nota} \labtt{mfq} \rm Let $M$ be a BW lattice of rank 
$2^{d-1} \ge 3$.   Let $Q$ be a lower group (see \ref{brw}) 
in
$Aut(M)$, i.e. in some $BRW^0(2^d,+)$ subgroup of $Aut(M)$, which
by induction is isomorphic to $BRW^0(2^d,+)$ or $d-1=3$ and
$Aut(M)\cong W_{E_8}$.    Also,let 
$f\in Q$ be a fourvolution, $F:=N_{Aut(M)}(Q) \cong BRW^0(2^d,+)$; see
\ref{brw}.   Now let $r$ be duality level of $M$ (see \ref
{dualitylevel}).  Then $r
\in
\{0,1\}$ and $r\equiv d (mod \ 2)$.  
\end{nota}

\begin{de}\labtt{nextbw} \rm {\it The Next BW Procedure.  } 
We use  notations $M, F, Q, f$ as in \ref{mfq}.   

Form $M_1 \perp M_2$, two orthogonal copies of $M$ based on
the isometries
$\psi_i:M \rightarrow M_i$ and let $V_i := \QQ \otimes M_i$ be their
ambient rational vector spaces.  Set $V:=V_1 \perp V_2$.  
Also,  we use
$\psi_{ij}:=
\psi_i^{-1}\psi_j$, the natural isometry from $M_i$ to $M_j$, extended to 
$V_i \rightarrow V_j$.   See Notation \ref{diagonalnotation}.

Define $Q_i, F_i$ and $f_i$ and the groups and element in $End(V_i)$
corresponding to $Q, F$ and $f$ under $\psi_i$.  Extend their actions
to $V$ in the natural way.      Also, define the group
$Q_{12}$ as the natural  diagonal subgroup of $Q_1
\times Q_2$ and element 
$f_{12}:=f_1f_2\in Q_{12}$ acting on $V$ (see \ref{diagonalnotation}).

For the
SSD sublattices $\mi {1-r}, \mij {1-r}, \mijp {-r}$, we denote the
respective SSD involutions by $t_i, t_{ij}, t_{ij'}$.  Observe that
$-1=t_1t_2=t_{12}t_{12'}$.  For convenience and symmetry, we define
$t_{i'}:=-t_i, t_{i'j}:=t_{ij'}, t_{i'j'}:=t_{ij}$.  Finally, we define 
$D:= \la t_1, t_2, t_{12}, t_{12'}  \ra \cong Dih_8$ and 
$R:=\la Q_{12}, D \ra\cong \explus d$.  So, $R=Q_{12}D$, a central
product.

Define $L_d :=\del$, and $R$-invariant lattice.  We call $L_d$ the {\it
\tbw -successor to 
$M$}.   From \ref {dottwists} and $M^*=M[-r]$, we deduce that $L_d$ is an 
integral lattice and since elements 
of the the above generating set have
even norms, $L_d$ is even.
\end{de}

\begin{lem}\labtt {dgbed} 
$\dg {L_d} \cong 1, 2^{2^{d-1}}$ as $d$ is even, odd.  Therefore, the
duality level is the remainder of $d+1$ modulo 2.  
\end{lem}
\pf
When $d$ is even, $L:=L_d$ is the kernel of the epimorphism 
$M_1 \perp M_2 \rightarrow M/M[1]$, defined by $(x^{\psi_1}, y^{\psi_2})
\mapsto x+y+M[1]$.  Since $M_1 \perp M_2$ is unimodular, $\dg L \cong
2^{2^{d-1}}$.  

When $d$ is odd, this is the same as the kernel of the epimorphism 
$M^*_1 \perp M^*_2 \rightarrow M^*/M$ defined by 
$(x^{\psi_1}, y^{\psi_2}) \mapsto x+y+M$.  
Since $M^*_1 \perp M^*_2$ has determinant 
$2^{-2^{d}}$, $L$ is unimodular.  

The statement about duality level follows from \ref{dualityleveltwist}.  
\eop

\begin{lem} \labtt {twistednextbw}  
We take 
$L := L_d$  
(in the notation \ref {nextbw}).  Here, $r \in \{0,1\}$, $d=rank(L)$
and $r \equiv d \ (mod \ 2)$.  Then: 

(i)  
$L$ is
the sum of any three of the four lattices 
$$ M_1[1-r],
M_2[1-r], M_{12}[-r], M_{12'}[-r].$$ 

(ii) 
$L^*$ is
the sum of any three of the four lattices 
$$M_1[1-r],
M_2[1-r], M_{12}[-1], M_{12'}[-1].$$  

\end{lem} 

\pf  This follows from \ref{dih8onlatt}.  Here is a different proof.  
(i) 
Let  $i=1$ or $2$.   
Define $j$ by: $\{1,2\}=\{i,j\}$.  

Now, we observe that for any integer $k$, 

(a) $\mi k \le \mj k + M_{12} [k] \le \mj  k + M_{12}[k-1]$;

(b) $M_{12'} [k] \le M_{12'} [k] + M_{12} [k] = M_{12} [k] + \mi {2+k}
\le M_{12} [k] + \mi {1+k}$.  

At once, (i) follows.  

For (ii), (a) and (b) prove equality of 
$N = M_1[1-r]\perp M_2 [1-r]+M_{12}[-1]$ and 
$N'= M_{12}[-1] + M_{12'}[-1] + M_i[1-r]$.  
It is clear by taking dot products that  $N=N'$ is in $L^*$. 
If $r=0$, $L=N$ and if $r=1$, $N/L \cong \monetwo {-1} / \monetwo 0 \cong
2^{2^{d-1}}$, whence $N=L^*$ (see \ref {dgbed}).  
\eop

\begin{coro}\labtt{density}
In the notation of \ref{commutatordense} and \ref{nextbw}, 
the $R$-module $L_d$ is commutator dense with respect to 
any fourvolution in $R$.   
\end{coro}
\pf
Since $Q_{12}$ acts diagonally on $V_i$, 
we deduce 
$[\mi k ,Q_{12}]=[\mi k,f] = \mi {k+1}$ 
for all $k$ and $i=1,2$.  
Also, 
$[\mij k , Q_{12}]=\mij {k+1}$.  
Note also, that $[\mi k , t]=\mij {k+1}$ and
$[\mij k , t ]= 0$ or $2 \mij  k = \mij {k+2}$ for 
$t=t_{ij}$ or $t_{ij'}$.  Similar statements hold for the
$\mijp k$.  
Since $R$ is generated
by
$Q_{12}$ and
$\la t_1, t_{ij}\ra$, 
$[L_d,Q_{12}]=[L_d,R]$.

We prove density 
 first prove for a few special cases of $f$.  

 Take $f=f_1f_2$, which acts diagonally.  
Then,  
by induction, 
$[\mij k , f_1f_2] = [\mij k , Q_{12}]=\mij {k+1}$ 
and 
$[\mi k, f]=\mi {k+1}$ and similarly for $\mj k$.  
So we have density for this $f$.

Now, let $f=t_1t_{12'}$.  
Then for $(x,y) \in V_1\perp V_2$,
$(x,y)(1-f)=(x-y,y+x)$.  
Since $(1-f)^2=-2f$, $L_d(1-f)$ 
contains $2L_d$
and the diagonal $\mij {1-r}$, which generate 
$[L_d,R]$ (see the first paragraph).  So, we have density for this $f$.

Finally, let $f$ be an arbitrary fourvolution in $R$. Then $|L:L(f-1)|$
has order $|L:2L|^\half$ so $L(f-1) \le [L,Q]$ implies that 
$L(f-1)=[L,Q]$.  
\eop

\begin{coro}\labtt{sultrycompatibility}  For $i=1,2$ and for all integers
$j \ge 0$, $M_i [1-r]\cap L[j]= M_i [1-r+j]$.
\end{coro}
\pf  Fix $i$.   
We choose $f_i$ for the twisting since it preserves the $M_i[k]$.  
The
equalities are valid for
$j$ even since $L[2k]=2^kL$, for all $k \ge 0$ and $M_i[1-r]$ is a direct
summand of $L$.  Now, 
$M_i [1-r]\cap L[1]\ge  M_i[2-r]$. Applying one more twist, which is a
scaled isometry, we get $M_i [2-r]\cap L[2]\ge  M_i[3-r]$.  Since 
$M_i [2-r]\cap L[2] = M_i[2-r] \cap 2L = M_i[2-r] \cap 2M_i[1-r] =
2M_i[1-r]=M_i[3-r]$, whence all our containments are equalities. 
\eop

We give a fairly complete account of minimal vectors.

\begin{lem}\labtt{minvec}  (i)  A minimal vector of $L_d$ has norm 
$2^{\lfloor {d\over 2}\rfloor}$ and is 
in $\mone{1-r}$ or $\mtwo{1-r}$ or has the form $x_1+x_2$, where each
$x_i$ projects to a minimal vector of $\mi{1-r}$, for $i=1,2$.   Its norm
is
$\mu (\mi{1-r})=2^{1-r}\mu (M) = 2^{\lfloor {d \over 2}\rfloor }$. 

(ii)  The minimal vectors span $L_d$.
\end{lem} 
\pf 
(i) 
Suppose that the minimal vector $x$ is not in $\mone{1-r}$ or
$\mtwo{1-r}$.  Write $x=x_1+x_2$, where $x_i$ is the projection of $x$ to
$V_i$, $i=1,2$.   Since $x_i \in \mi{-r}$, $x_i$ has norm at least $\half
\mu (\mi{1-r})$, whence $(x,x) \ge \mu (\mi{1-r})$.  It follows that these
inequalities are equalities.  The last statement follows from induction
and
\ref{twistednextbw}.  

Easily, (i) implies 
(ii) since $L_d$ is
the sum of three sublattices spanned by minimal vectors.  
\eop  

\begin{coro}\labtt{ldisbwt} $L_d$ is a lattice of BW-type. 
\end{coro}

\begin{de} \labttr{type123} A minimal vector in $L_d$ has {\it type 1, 2
or 3}, respectively, as it is in $M_1[1-r], M_2[1-r]$ or in neither.  
These three types partition $\mv {L_d}$.  
\end{de}

\begin{lem} \labtt{beindec} For $d \ge 2$, 
$L_d$ is an indecomposable
lattice. 
\end{lem}
\pf
As in the proof of \ref{kstlattices}, we see that the minimal vectors of
$L$ are partitioned into equivalence classes by membership in the $L_i$. 
However, it is clear from \ref{minvec} and induction that there is just
one equivalence class in the sense of \ref{kstlattices}. 
\eop

The following terminology will be useful.  It applies to lattices used in
\ref{nextbw} and later.  

\begin{de}\labttr{ssbw} 
A lattice $M$ is a {\it scaled BW-lattice}, abbreviated {\it sBW
lattice},  if there is an integer
$s>0$ so that $M\cong \sqrt s \bw e$, for some $e>0$.  
 A sublattice $M$ of a BW-lattice $L$ is called a
{\it suitably scaled Barnes-Wall sublattice (relative to $L$)},
abbreviated  {\it ssBW sublattice}, if  $M$ is a sBW lattice and $\mu
(M)=\mu (L)$.  

We use the notation $\bw {p,q}$, $p \le q$, for a scaled copy of $\bw p$
whose isometry type is suitable as a sublattice of $\bw q$, i.e. 
a sBW lattice with minimum norm $2^{\lfloor {q\over 2}\rfloor }$.  
\end{de}

\section{The groups $\rd d$, $\gd d$ and invariant lattices}

\begin{nota}\labttr{ldd} 
In this section, $L=L_d$  and $d\ge 2$ have the meaning of \ref{mfq} and
\ref{nextbw}.
\end{nota} 

\begin{de}\labtt{rdgd}\rm 
We define $R:=\rd d:=\la Q_{12}, t_i, t_{ij} \ra \cong \explus d$, 
where
$Q_{12}$ and the involutions are as in 
\ref{nextbw}. 
We define $\gd d:=N_{Aut(L)}(\rd d)$.  
\end{de}

\begin{de}\labtt{lowerstuff}\rm 
Elements and subsets  of $\gd d$ are called {\it lower} if in $R$
and are otherwise called {\it upper}.  
In particular, a fourvolution \ref{fourvolution} may be called upper
or lower.  

\end{de}

\begin{thm}\labtt {shapegd}For $d\ge 2$, 
$\gd d \cong  BRW^0(2^d,+) \cong \ratholoex d$.  
\end{thm} 
\pf  The cases $d\le 3$ have been discussed earlier 
(and the case $d=4$
was treated explicitly in \cite{POE}).  We may assume that $d\ge 4$.   
Since $\gd d$ is finite, containing
$R$ as a normal subgroup, 
$\gd d$ is contained in
$\widetilde G = BRW^0(2^d,+)$, the natural
$\ratholoex {d}$ subgroup of $GL(2^d,\CC )$ containing $R$;  
see Appendix A2.)

Let $D$ be the dihedral group of order 8 described in
\ref{nextbw}.  Then $D \le G:=Aut(L)$.   Let $t \in D$ be a
noncentral involution.   We claim that
$C_{\gd {d}}(t)R/R$ corresponds  to a maximal parabolic in $\widetilde
G/R$. For standard theory about parabolic subgroups, see \cite{Car}.

Suppose $t=t_1$ or $t_2$.  
By induction, $Aut(\mi {1-r})$ contains a copy of $\gd {d-1}$ as
$N_{Aut(\mi {1-r})}(Q_i)$,  and $S:=Stab_G(\mone{1-r} \perp
\mtwo{1-r})$  contains a group $T$ of the form $[\explus {(d-1)}
\times
\explus {(d-1)}].[\Omega^+(2(d-1) \times 2]$.  
Also, since $t_{ij}$ interchanges $\mone {1-r}$ and $\mtwo {1-r}$, 
it normalizes 
this group.  Its image in  $\widetilde G/R$ is a maximal parabolic, the
stabilizer of a singular vector.

Suppose that $t=t_{ij}$ or $t_{ij'}$.  Then the above argument goes
through with $\mi{1-r}$, $t_i$ replaced by $\mij {-r}$, $t_{ij}$, and
gives a distinct subgroup of the form 
$[\explus {(d-1)}
\times
\explus {(d-1)}].[\Omega^+(2(d-1), 2)\times 2]$ containing $R$.  
(Proof of distinctness: in both cases, the center of the respective
stabilizer is $\{\pm 1, \pm t \}$.)

Therefore,  
$G/R$ contains two different maximal parabolics of $\widetilde G/R$, whence
$G=\widetilde G$, \cite{Car},  so we are done.  
\eop

\begin{rem}\labtt {avoidog}\rm  Note that the \ref{shapegd} uses 
only a basic result about 
orthogonal groups (maximality of certain stabilizers) but nothing very
explicit about their interior structure, nor about particular elements. 
This is possible since we have a suitable uniqueness statement. 
\end{rem}

\begin{lem}\labtt {newlowergroup}
The subgroup of $Aut(L)$ which is trivial on $L/L[1]$ is just $R$.   In 
the notation \ref {nextbw}, $L[1]=\mone {-r} + \mtwo {-r} + \monetwo
{1-r}$.  
\end{lem}
\pf 
Let $T$ be the subgroup trivial on $L/L[1]$.  Note that $[L,R]=L[1]$,
\ref{twistednextbw}.  
Therefore, 
$T \ge R$.  

Assuming $T > R$, 
we have a normal nontrivial 2-group $T/R$ in $\gd d/\rd d$.  Since the
latter quotient is simple, the shape of $\gd d$ given in
\ref{shapegd} shows that this is impossible.  
\eop

\begin{nota}\labttr {sf}
For $x \in \mv L$, let $SF(x) :=x^R$.  Call this the 
{\it sultry
 frame
containing
$x$. }(See \ref{sultryinterpretation}).   From
\ref{newlowergroup} and the structure of 
$R\cong \explus d$, $SF(x)$ is a
double orthogonal basis, of cardinality $2^{d+1}$.   
\end{nota}  

\begin{prop} \labtt{sfequivalence}  Let $x, y \in \mv L$.  Equivalent
are (i) $y \in SF(x)$; (ii) $x-y \in L[1]$.
\end{prop}
\pf  Trivially, (i) implies (ii).  For (ii), we use 
 a familiar argument.  Let $z \in SF(x)$.  
First we note that $z\pm y \in L[1]$ is 0 or has norm at least $2\mu (L)$.
Assuming $y \ne \pm z$, we have $(z\pm y, z\pm y)=2\mu(L)\pm 2(z,y)\ge
2\mu (L)$, whence $(z,y)=0$.  This is not the case for every $z \in SF(x)$.
\eop

\begin{prop}\labtt{numberminvect} The number of minimal vectors is
$(2^d+2)(2^{d-1}+2)\dots (2^2+2)(2+2)$.  The values for small $d$ are:

$$\begin{matrix}     d&|\mv L | & \hbox{ Prime Factorization} \cr
0&2&2\cr
1&4 &   2^2  \cr2&24&2^3 3\cr  
3&240  &  2^43.5  \cr
4&4320& 2^53^35
      \cr
5&146880& 2^6 3^3 5. 17
      \cr
6&9694080& 2^73^45.11.17 
     \cr
7&1260230400   & 2^83^45^211.13.17  
\cr
8&325139443200& 2^93^55^2.11.13.17.43\cr
9& 167121673804800&
2^{10} 3^5 5^2  11.13.17.43.257  \cr 
10 &171466837323724800 & 
2^{11}3^8 5^2 11.13.17.19.43.257 \cr
11& 351507016513635840000 & 
2^{12}3^8 5^4 11.13.17.19.41.43.257\cr 
12& \ 1440475753672879672320000\ \  &
2^{13}3^9  5^4 11.13.17.19.41.43.257.683\cr 
\end{matrix}$$
\end{prop}
\pf
Use \ref{minvec}, \ref{sf}, \ref{sfequivalence} and induction.  
\eop

\begin{coro} \labtt{sfstab} If $x\in \mv L$, $Stab_{\gd d}(x+L[1])/\rd d$ is a
maximal parabolic subgroup of $\gd d/\rd d$ of the shape 
$2^{d \choose 2}{:}GL(d,2)$.
\end{coro}
\pf
The pairs $\{\pm x\}$ of minimal vectors in this coset is an orbit of $\rd d$
for which a point stabilizer $E$ is elementary abelian of order
$2^{1+d}$; see \ref{sfequivalence}.  These pairs of vectors are
exactly the minimal vectors of the total eigenlattice of $E$, so as a
set are stable under
$N_{\gd d}(E)$, which has the indicated properties.  
\eop

\begin{coro}\labtt {o2deven} When $d$ is even, 
$L \cap 2L^*=L[1]$, whence the lower group is normal in $Aut(L)$
and $\gd d=Aut(L)$.  
\end{coro}
\pf 
When $d$ is even, the duality level is 1, whence $L \cap 2\dual
L=L \cap 2L[-1]=L \cap L[1]=L[1]$ is invariant by the entire automorphism
group.  Now use \ref{newlowergroup}.
\eop

When $d$ is even, this result essentially solves the problem of
determining the automorphism group. 
 The case $d$ arbitrary is harder (it
is finally proved in
\ref{g=gd}, which does not use
\ref{o2deven}).

\begin{prop}\labtt {reductionmodp}
For an odd prime, $p$, $L/pL$ is an absolutely irreducible module for
$\gd d$.  
\end{prop}
\pf  This is trivial, since $R$ acts absolutely irreducibly.  
\eop

\begin{lem}\labtt {irreducibility} 
For all $k$, $L[k]/L[k+1]$ is an absolutely irreducible $\FF_2$-module
for 
$\gd d$.  
\end{lem}
\pf    This is easy to check for $d \le 4$, so we assume that $d \ge 5$
and use induction.    

We may assume that $k=0$.   
Let $D \cong Dih_8$ be as in \ref{nextbw}.  
For a noncentral involution
$t$ of $D$, we get by induction that $C_{\gd d}(t)$ acts irreducibly
on each 
$L^\pm (t) / [ L^\pm (t), C_R(t) ]$.

Since $L$ is a sum of fixed point sublattices for the 
noncentral involutions of $D$ \ref{equivcd}, it follows that $L/L[1]$
has two absolutely  irreducible composition factors for $C_{\gd
d}(t)$, each of dimension
$2^{d-2}$.   

The group $C_{\gd d}(t)$, of shape 
$[\explus {(d-1)}
\times
\explus {(d-1)}].\Omega^+(2(d-1), 2)$ (discussed in the proof of
\ref{shapegd}), acts on $L/L[1]$ and has exactly one irreducible
submodule, namely $Tel(L,t)/L[1]$, and two composition factors
(this follows by induction on $d$, since on any eigenspace for $t$,
we know the irreducible quotients for any lattice invariant under
$C_{\gd d}(t)$).   These irreducible submodules distinct as
$t$ ranges over a set of generators for $D$ (e.g. $t_1, t_{12}$ for the
group of 
\ref{nextbw}).  

It follows that $L/L[1]$ is irreducible for the action of $\gd d$. We
now prove absolute irreducibility.
If   $K$
is an extension field of $\FF_2$ and  $K\otimes L/L[1]$ decomposes,
then its restriction to $C_{\gd d}(t)$ would have over 4 composition
factors (since $O_2(C_{\gd d}(t))$ acts nontrivially),  which is
impossible  since, by induction, the composition factors for
$C_{\gd d}(t)$ are absolutely irreducible of dimension $2^{d-2}$.    
\eop

\begin{prop} \labtt {allinvariantlattices}
Let $M$ be a lattice in $\QQ \otimes L$ which is invariant under $\gd d$. 
Then there is a rational number $r$ so that $rM=L$ or $L[1]$.  
\end{prop}
\pf
We may assume that $M \le L$.  By \ref{reductionmodp}, we may assume that
$L/M$ is a power of 2.   For some positive integer $n$, $2^nL \le M$.  Now
use \ref{irreducibility} and the fact that $[L[k],\rd d]=L[k+1]$.   
\eop

Later, in \ref{g=gd}, we prove that $\gd d$ is all of $Aut(L_d)$.

\section{Sultriness}

When $f$ is a fourvolution on a lattice $L$, $1-f$ 
(actually, any of $\pm 1 \pm f$) is an endomorphism of $L$  
which is also an isometry scaled by $\sqrt 2$.
 Next, we see that a sulty
transformation is naturally interpreted as a {\it scaled lift of a
transvection}, a point which suggested the term ``sultry''.   

\begin{thm}\labtt{sultryinterpretation} 
The function $\gamma_f : Aut(L) \rightarrow Aut(L)$, 
$x \mapsto (1-f)^{-1}x(1-f)$, normalizes $R=\rd d$ and $\gd d$.
Furthermore, $\gamma_f$ is the identity on $C_R(f)$ and if 
$x \in R\setminus C_R(f)$, 
$\gamma_f$ takes $x$ to $fx \in R\setminus C_R(f)$, hence normalizes
and induces an outer automorphism on the dihedral group $\la f, x \ra$. 
Hence, on $R/Z(R)$, $\gamma_f$ acts as the transvection associated to the
nonsingular point $Z(R)f$ of $R/Z(R)$.  
\end{thm}
\pf
Since ${1 \over \sqrt 2}(\pm 1\pm f)$ is orthogonal, the image of
$\gamma := \gamma_f$ is a subgroup of the orthogonal group.  Since any
$\pm 1\pm f$ carries each $L[k]$ onto $L[k+1]$, the image of $\gamma$
stabilizes $L$.  We conclude that $\gamma$ takes $Aut(L)$ onto itself.  

We calculate that $x(1-f)=x-xf=x-f^{-1}x=x+fx=(-f+1)(fx)$, which
proves the remaining statement.  
\eop

\section{Proof of uniqueness  }

\def\ldg#1{{\bf (LDG #1)}  }  
\def\ldl#1{{\bf (LDL #1)}  } 
\def\xol#1{{\bf (XOL #1)}  } 
\def\hypxd{{Hypothesis Xd}} 
\def\hypxdmo{{Hypothesis $X(2^{d-1})$}} 
\def\hypxdmt{{Hypothesis $X(2^{d-2})$}} 
\def\hypxdpo{{Hypothesis $X(2^{d+1})$}} 
\def\hypxdpt{{Hypothesis $X(2^{d+2})$}}

\def\xd#1{{\bf (Xd.#1)} } 
\def\x#1#2{{\bf (X#1.#2)} }


\begin{nota}\labttr {frakx} 
Given $d\ge 3$ and  $L_1, L_2$, we let 
${\frak X} := {\frak  X}(L_1,L_2)$ be the set of all
$X$-quadruples of the form $(L,L_1,L_2,t)$;  see \ref{condx}.  
\end{nota}

\begin{thm} \labtt{bwtclassification} 
We use the notation in \ref{condx}, \ref{mfq}, 
\ref{nextbw},
\ref{rdgd} and 
\ref{frakx}. 
Suppose that $d \ge 3$ and $(L_1, L_2)$ is an orthogonal pair of
lattices, so that each 
$L_i$ is BW-type of rank $2^{d-1}$.  

(i)  ${\frak X}$ is an 
orbit under the natural action of $F_1 \times F_2$, where
$F_i:=Stab_{Aut(L_i)}(L_i[1-r])$ (see \ref{nextbw}; by
induction,
$F_i
\cong \gd {d-1}$).  Define $Q_i:=C_{F_i}(L_i/L_i[1])$.  

The elements of $\frak X$ are in correspondence with each of the following 
sets.

(a) $F_1/Q_1$;  

(b) $F_2/Q_2$;

(c) Pairs of involutions $\{s, -s\}$ in the orthogonal group on $V$ which
interchange
$L_1$ and $L_2$.

(d) Dihedral groups of order $8$ which are generated by the SSD
involutions associated to $L_1, L_2$ and involutions as in (c).

(ii) (a) The subgroup $G_L^0$ of $F_1 \times F_2$  which stabilizes
$L$ has structure 
$Q_1 \times Q_2 \le G_L^0$ and $G_L^0/Q_{12}$ is the diagonal subgroup of 
$F_1/Q_1 \times F_2/Q_2$ with respect to the isomorphism induced by
$s$, an involution as in (i.c).

(b) 
The subgroup $G_L$ of $Aut(L_1 \perp L_2)
\cong Aut(L_i)
\wr 2$ which stabilizes $L$ is
$G_L^0\la  s \ra$.  
We have $G_L \cong [\explus {d-1} \times \explus
{d-1}].[\Omega^+(2(d-1),2)
\times 2]$.  

(c) The subgroup of $G_L$ which acts trivially on $L/L[1]$ is $R:=\la
Q_{12}, s, t_i \ra$, where $t_i$ is the SSD involution associated to
$L_i$.  The quotient
$G_L/R
\cong 2^{2d-2}{:}\Omega^+(2d-2,2)$
is a maximal parabolic subgroup of $Out^0(\explus {d})
\cong
\Omega^+ (2d,2)$. (See Appendix A.0).  
\end{thm}

The extension in 
(c) is split, despite
$\gd e$ being nonsplit over $\rd e$ for $e \ge 4$.  See Appendix A2.

\pf  (i) 
We prove the classification by induction.  For $d=2$,
$Aut(L_{D_4})\cong  2^{1+4}_+[Sym_3 \wr 2]$ and for $d=3$, 
$Aut(\leh )\cong
W_{E_8}$.  
When $d=4$, the main theorem follows from the arguments of \cite{POE}.  

For the rest of the proof, we assume that $d \ge 4$.  
By induction, a lattice satisfying the $X(2^{d-1})$ condition is uniquely
determined up to isometry.  
This applies to the lattices $L_1,  L_2$.

Let $L$ be any member of $\frak X$ and set $G:=Aut(L)$.  There is an 
X-quadruple 
$(L,L_1,L_2,t)$.    
Then $det(L)$ and $|L{:}\mone {1-r} \perp \mtwo
{1-r}|$ are determined.

Let $p_i$
be the orthogonal projection of
$L$ to $V_i:=\QQ \otimes L_i$, for $i=1,2$.  
Then $L^{p_i}$ is a lattice containing 
$L_i$ with quotient isomorphic to  $2^{2^{d-2}}$.
Since $Q$ acts trivially on $L/[L_1\perp L_2]$, $Q_i$ acts trivially on
$L^{p_i}/L_i$.  Therefore, $L^{p_i}$ is the $-1$ twist of $L_i$ with
respect to $Q_i$, i.e., $L^{p_i}=L_i(1-f_i)^{-1}$, for a suitable
fourvolution $f_i \in Q_i$.

There is a dihedral  subgroup $D$  of $R$ so that $t \in D$. 
If $y \in D$ is an involution which 
does not commute with $t$, then 
$y$ interchanges $L_1$ and
$L_2$.  Also, if $J^\pm$ are the eigenlattices for $y$, then 
$L$ is the sum of any three of the four $L_1, L_2, J^+, J^-$, by
\ref{dih8onlatt}.

It follows that $L$ is determined by $D$ in the sense that 
$L=[L_1 \perp L_2] + ([L_1 \perp L_2]^+(y))(f-1)^{-1}$, where
$f=ty$.

is
the sum of the fixed point sublattices of the involutions of $D$ (see
\ref{dih8onlatt}, 
\ref{twistednextbw}).

Now, to what extent does 
$L_1 \perp L_2$ determine $D$?  The answer is: up to conjugacy in 
$Aut(L_1 \perp L_2) \cong \gd {d-1} \wr 2$ (note that $d\ge 2$ here).  
Our group $D$ is generated by the center of the natural index 2 subgroup
of $Aut(L_1 \perp L_2)$ and a wreathing involution.  
In general, 
wreathing involutions in a wreath product of groups $K \wr 2$
form an orbit under the action of either direct factor isomorphic to $K$
in the base group of the wreath product. This proves correspondence with
(c) and (d).  The stabilizer subgroup is diagonal in the base group $K
\times K$, and either direct factor represents all cosets of the stabilizer
(whence the equivalence with (a) and (b))

It follows that, up to isometry preserving $L_1 \perp L_2$,
$D$, hence $L$, is determined by the pair of indecomposable lattices
$L_1$ and
$L_2$.

Proof of statements (ii) and (iii) are easy.  The statement about
parabolic subgroups is proven with a standard result from  
the theory of
Chevalley groups, e.g.
\cite {Car}.  Independently of that theory, the maximality could be proved
directly by showing that there is no system of imprimitivity on the set
of isotropic points.    This is an exercise with Witt's theorem. 
\eop

\section{Minimal vectors, the zoop2 property and $Aut(\bw d )$}

We continue to use the notations of \ref{ldd} and \ref{rdgd}.

\begin{rem}\labtt {twistminimalvectors}\rm 
For $L=\bw d$ a Barnes-Wall lattice and $k \in \ZZ$, we have $\mv {L[k]} =
\mv L [k]$ (see
\ref{dottwists}).  
\end{rem}

\begin{thm} \labtt{onminvect} We use notation of \ref{nextbw}.   
The group $\gd d$ acts transitively on the set of minimal vectors.

\end{thm}
\pf   We use notation of \ref{type123}.  
It is clear that the minimal
vectors of types 1 and 2 are in a single
$\gd d$-orbit, say $\cal O$.  Consider a minimal vector
$u+v$ of type 3.  
We assume that $L$ corresponds to the involution $s=t_{ij'}$ in the
sense of \ref{bwtclassification} (i). 
Then $v$ and $u^{t_{i'j}}$
differ by an element of
$\mtwo {1-r}$, so by induction, these are in the same orbit under
$Q_2$, equivalently under $Q_{12}$.   Therefore,  $u$ and $v$ are in
the same 
$R$-orbit.

By induction, we have transitivity on the minimal vectors of type 3
by the group $RF_{12}$, where the second factor is the natural diagonal
subgroup of $F_1 \times F_2$.  Call ${\cal O}'$ the
orbit containing the type 3 minimal vectors.  

Suppose that $\gd d$ is not transitive.  {\it First Contradiction.} 
Then $\mv
L$ is the disjoint union of two orbits 
$\cal O$ and
${\cal O}'$ and so $\gd d$ preserves the $\ZZ$-span of $\cal O$, which is just 
$\mone {1-r} \perp \mtwo {1-r}$, an orthogonal sum of two orthogonally
indecomposable lattices.  Thus $\gd d' \ge R$ leaves both summands
invariant, which is impossible since $R$ is irreducible on $\CC \otimes
L$.  
{\it Second Contradiction. }  The lower involutions form a conjugacy class
in
$\gd d$, so there is $g \in \gd d$ which conjugates $t_1$
 to $t_{12}$.  Then $g$ takes the set of minimal vectors fixed by $t_1$
(those of type 2) to those fixed by $t_{12}$, which are contained in those
of type 3.  Therefore $\cal O$ and ${\cal O}'$ are not distinct
orbits.  Transitivity follows.  
\eop

Now we give a few results about stablilzers in $\gd d$.  These will be
strengthened later.

\begin{lem}\labtt {diagonallowergroup}
(i) If $F$ is a sultry 
frame and $x \in F$, then $\{ g\in R|x^g=\pm x\}=\{g
\in R| y^g=\pm y \hbox{ for all }y\in F\}$ is a maximal elementary abelian
subgroup of $R$.  Call it $R_F$.  The quotient $R/R_F$ operates
regularly on the eigenlattices.

(ii) Define $C_F:=C_{\gd d}(F/\{\pm 1\}) := \{ g \in \gd d | y^g=\pm y \hbox{
for all }y \in F\}$.  This is elementary abelian and  has shape
$2^{1+d+{d
\choose 2}}$.

(iii) Its normalizer $N_F:=N_{\gd d}(C_F)=Stab_{\gd d}(F)$ in $\gd d$
satisfies
$N_F/C_F
\cong AGL (d,2)$.  We have $R_F\le C_F$.  
\end{lem}  

\pf  (i) Set $P:=\{ g\in R|x^g=\pm x\}, Q:=\{g
\in R| y^g=\pm y \hbox{ for all }y\in F\}$.  Observe that $Q\le P$
and $Q$ is elementary abelian.  Transitivity of $R$ on $F$ and
normality of $P$ in $R$ implies that $P=Q$ has order $2^{d+1}$.

(ii) This follows from (i) and the order of the unipotent radical for the
stabilizer of a maximal totally isotropic subspace for $\Omega^+(2d,2)$.

(iii) This follows from the actions of $R$ on $R_F$ together with the
structure of the 
stabilizer of a maximal totally isotropic subspace for $\Omega^+(2d,2)$. 
\eop

\begin{de}\labttr{zoop2} Suppose that $F$ is a frame.  A subset $S
\subseteq V$ has the  {\it zop2 property (with respect to $F$)} 
if $|(x,y)|$ is 0 or a power of 2, for all $x \in F$ and
$y
\in S$.   We say that $S$  has the {\it zoop2 property } if is has the
zop2 property and just one  power of 2 occurs among the scalars $|(x,y)|$,
$x\in F, y \in S$.  
\end{de}

\begin{lem}\labtt{minvectzoop2}   For all integers $p, q$, any minimal
vector of
$L[q]$ has the zoop2 property with respect to any sultry frame in
$L[p]$.  
\end{lem}  
\pf 
We may assume that $p=0$ and $q\in \{0, 1\}$.  The property is easy to
check for
$d\le 3$.  We assume $d\ge 4$.   Say $x \in \mv L$,
$y \in SF(x)$ and
$z \in \mv {L[q]}$ so that $(x,z)\ne 0\ne (y,z)$.  Take a lower involution
$t$ so that
$t$ fixes
$x$ and
$y$.  Then $x, y \in L^+(t)$, a ssBW lattice, and $z$ projects to a
minimal vector in $L^+(t)[q]$, so we are done by induction.   
\eop

\begin{de}\labttr{support} Given a sultry frame $F$ of $2^{d+1}$ elements,
there are $2^{d}$ subsets which form a basis.  Suppose $X$ is such a set. 
If
$v\in V$, we write
$v=\sum_{x \in X} a_x x$ and define the {\it support of $v$} to be the set 
$\{x\in X| a_x \ne 0 \}$.  This depends on the double basis  $F$, not on
the choice
$X
\subset F$.
\end{de}

\begin{lem}\labtt {spanninghamlatt}  
Suppose that $d \ge 2$.  
Let $x \in \mv L$ and $A(x):=\{ y \in \mv L | (x,y)=\half (x,x) \}$.  
Then $A(x) \cup \{x\}$ spans a lattice isometric to the Hamming code
lattice 
described in 
\ref{lbhame}.  In particular, there is a labeling of $SF(x)/\{\pm 1\}$ with
$\FF_2^d$ so that the elements of $A(x)$ have support which is an affine
2-space.  
\end{lem} 
\pf
Define $J_0$ to be the square lattice spanned by $SF(x)$ and $J$ the
lattice spanned by $A(x)$ and $SF(x)$.  In a natural way, $J/J_0$
corresponds to a nonzero code in $\half J_0/J_0 \cong \FF_2^{2^d}$.  
Since
$\mu (L)=(x,x)$, this code has minimum weight at least 4. 
For $y \in A(x)$, $supp(y)$ is a 4-set with respect to the double basis
$SF(x)$, \ref{minvectzoop2}.  Therefore the minimum weight of $C$ is
4.  

Note that we have an action of $AGL(d,2)$ on $\half J_0/J_0 $ 
by coordinate permutations.  This follows from \ref{diagonallowergroup}.
This action is triply transitive.  Since $C$ has minimum weight 4, its
weight 4 codewords forms a Steiner system with parameters $(3,4,2^d)$
which is stable under this action of
$AGL(d,2)$.  Such a system is unique since in $AGL(d,2)$, the stabilizer
of three points fixes a unique fourth point.  
Therefore, $C$ is the code $\hame d$, up to equivalence.

Finally, we must show that $A(x) \cup \{x\}$ spans $J$. 
If $d=2$, $L \cong L_{D_4}$ and the result is easy to check directly. 
We assume $d\ge 3$.  
Let
$y
\in SF(x), y\ne \pm x$.  Since $d \ge 3$, we may choose a lower
involution
$t$ which fixes both
$x$ and $y$ (in the notation of \ref{diagonallowergroup}, $t \in R_F$). 
Let $L^+$ be the sublattice of points of $L$ fixed by $t$, a sBW
lattice.   Then, induction implies that the sublattice of $L^+$ spanned by
$A(x)
\cap L^+$ contains $y$.  We conclude that $SF(x) \subset span(A(x) \cup
\{x\})$, and we are done.  
\eop

\begin{prop} \labtt{stabminvec} 
Suppose that $d \ge 4$.  
For $x \in \mv L$, 
$$Stab_{Aut(L)}(x) \le Stab_{Aut(L)}(SF(x)).$$   
\end{prop}
\pf
Define $A(x):=\{ y \in \mv L | (x,y) =\half (x,x)  \}$.  
By \ref{spanninghamlatt}  lattice $J:= span(A(x) \cup \{ x\})$
contains 
$SF(x)$ and is a copy of the lattice in \ref{lbhame}.
Since $d\ge 4$, 
given a weight 4 codeword, there exists another weight 4 codeword
which meets it in a 1-set (this is not so for $d=3$). 
Therefore, $SF(x)=\{ z
\in
\mv J | (z,J) 
\le
\half (z,z) \ZZ
\}$, we are done (see the proof of \ref{lbhame}).  
It follows from \ref{kstlattices} that  $Stab_{Aut(L)}(x) \le
Stab_{Aut(L)}(J)$.
\eop

\begin{coro}\labtt {g=gd}  For $d \ge 4$, 
$Aut(L_d)=\gd d$.
\end{coro}

\pf This follows since $\gd d$ is transitive on minimal vectors and the
stabilizer of some minimal vector in $G$ is contained in $\gd d$.
\eop

\begin{de}\labtt{layers}\rm Let $x\in \mv L$ and $SF(x)$ its sultry
frame.  Let $q\in \ZZ$ and $k \in \ZZ$.   
Define $A(L,x,q,k):=\{ z \in \mv
{L[q]} | (z,y)\in
\{ 0,
\pm 2^k\}\hbox{ for all }y \in SF(x) \}$.  This is is the {\it level $k$
layer in $\mv {L[q]}$ with respect to $x$} or $SF(x)$.  
\end{de}

\begin{lem}\labtt{invariantsection}  Suppose that the group $G_0$
factorizes as 
$G_0:=G Z$, where $G, Z$ are subgroups so that 
$[G,Z]=1$.  Suppose that $G_0$  acts on the set $X$ and that  $G$
stabilizes and acts transitively on a set of
$Z$-orbit representatives.   Let $S$ be
the set  of all $G$-invariant sets of 
orbit representatives.  Then $Z$ acts
transitively on $S$.  
\end{lem}
\pf 
A member of $S$ is determined by any element of $X$ which it  contains. 
Therefore the members of $S$ partition  $X$.  

Suppose that $X_1, X_2 \in S$.  Let $A$ be a $Z$-orbit and let $a_i$ be
the unique element of $A\cap X_i$, for $i=1,2$.  Take $z\in Z$ so that 
$a_1^z=a_2$.  Then $X_2$ and $X_1^z$ are both $G$-invariant sets
of orbit representatives and contain $a_2$, hence are equal.
\eop

Note that the next result deals with minimal vectors in all sultry twists
of $L$.

\begin{nota}\labttr{ioexp} For integers $d\ge 2$ and $p,q \in \ZZ$,
define $I(d,p, q)$ to be the set of integers  listed below. 
Here, $r\in \{0,1\}$ is the remainder of $d$ modulo 2, $m:=\lfloor
{d\over 2}\rfloor$ and $s\in \{0,1\}$ is the remainder of $p-q$ modulo
$2$.  We define

$$I(d,p,q):= \lfloor {p+1\over 2}\rfloor + 
\lfloor {q+1 \over 2} \rfloor  + \{-rs,0,1,\dots , m\}.$$ 
As usual $a+\{b,c , \dots \}$, means $\{a+b, a+c, \dots \}$.  We call 
$I(d,p,q)$ the {\it interval of exponents for dot products of 
minimal vectors}.  (See \ref{checkioexp} for an explanation of
this term.)   
\end{nota}

\begin{ex}\labttr{exioesp} 
Some  examples: $$I(3,p,p)= \begin{cases}  \{p,p+1\} & p \hbox{
even};\cr  
\{p+1,p+2\} & p \hbox{ odd.} \end{cases};  
I(3,p,p+1)=\{p,p+1,p+2\}; $$  
$$I(4,p,p)=\begin{cases} \{p,p+1,p+2\} & p \hbox{ even}; \cr 
 \{p+1,p+2,p+3\} & p \hbox{ odd.}\end{cases} ; 
I(4,p,p+1)=\{p+1,p+2,p+3\}.$$ 
\end{ex}

\begin{lem}\labtt{checkioexp}  The set of integers 
$$\{ (x,y)|x\in \mv {L[p]}, y \in \mv {L[q]} \}$$ is $\{ 0, \pm
2^k|k \in I(d,p,q) \}$; see \ref{ioexp}.    
\end{lem}
\pf  Let $m:=\lfloor {d \over 2}\rfloor$ and $r:=d-2m\in \{0,1\}$. 
  
First we take $p=q=0$.  Then for $x\in \mv L$,  with respect to a
basis $\Omega \subseteq SF(x)$ an element $y \in \mv {L[q]}$ has the form 
 $y=2^{-t} \sum_{u \in A} u\varepsilon_B$, where $A \subseteq \Omega$ is an
affine subpace of
$\Omega$ of dimension $a \le d$, and $a$ satisfies $m=m-2t+a$, or
$a=2t$.    Then
$(x,y)=0$ or
$\pm 2^{m-t}$.  
Thus, the values of  $a$, $t$, $m-t$  which occur are 
$\{0,2,\dots 2m \}$, 
$\{0,1,\dots ,m\}$, 
$\{0,1,\dots, m\}$,  
respectively.  

Next, if $p=0$ and $q=-1$, a similar discussion applies, but here $\mu (L[-1])=\half
\mu (L)$, so we get
the condition 
$m-1=m-2t+a$, or $a=2t-1$,  whence odd parity for $a$.   
Therefore, the   values of 
$a$, $t$, $m-t$  
 which occur are  
$\{1, 3, \dots , d+r-1\}$, 
$\{ 1,\dots, d-m\}$, 
$\{2m-d, 2m-d+1, \dots , m-1 \}$,
respectively.  

Suppose that $p=-1$ and $q=0$.  We get the condition 
$m=m-1-2t+a$, or $a=2t+1$ is odd.  Therefore the values of $a, t, m-t$
which occur are 
$\{1, 3, \dots , d+r-1\}$, 
$\{0, 1,\dots, d-m-1\}$, 
$\{2m-d+1, 2m-d+2, \dots , m\}$,
respectively.

For the general case, just observe that 
$I(d,p,q+2)=1+I(d,p,q)$, 
 $I(d,p+2,q)=1+I(d,p,q)$, and $I(d,p+1,q+1)=	1+I(d,p,q)$.    
\eop

\begin{lem}\labtt{cohomology} 
Let $G=AGL(d,2)$ act naturally on the permutation module
$A:=\FF_2^\Omega$, where $\Omega:=\FF_2^d$ with the
natural $G$-action.  Let $B$ be the submodule generated by the affine
subspaces of codimension 1.  For $d\ge 3$,  $H^1(G,B)=0$.  
\end{lem} 
\pf 
We have an exact sequence $0 \rightarrow B \rightarrow A \rightarrow A/B
\rightarrow 0$.  From this, the long exact cohomology sequence gives 
the exact sequence $H^0(G, A/B) \rightarrow H^1(G, B) \rightarrow H^1(G,
A)$.  The right term is, by the Eckmann-Shapiro lemma, isomorphic to 
$H^1(G_0, \FF_2 )$, where $G_0 \cong GL(d,2)$ is the stabilizer of 0 in
$G$. This is isomorphic to $Hom(G_0, \FF_2)$, which is trivial for $d\ge
3$.   The module $A/B$ is indecomposable for $AGL(d,2)$, with a faithful
module of dimension $d$ as the socle and quotient the trivial 1-dimensional
module.  Since the fixed points are 0, $H^0(G, A/B)=0$.  Exactness
implies that   $H^1(G,B)=0$.  
\eop

\begin{thm} \labtt{labeling}  Let $d\ge 2$ and let $F:=SF(x)$, for a
minimal vector
$x
\in L[p]$ (see \ref{sf}).  Also, define $H:=Stab_{\gd d}(F)$ (denoted
$N_F$ in \ref{diagonallowergroup} (ii) ).  

(i) There exists a 
basis $X$ contained in $F$ and labeling of $X$ by $\FF_2^d$ so
that with respect to $X$, $H$ is the 
monomial group  ${\cal E}_{{\cal C}_X}{:}AGL(d,2)$ (see
\ref{sumsandsigns}), where 
${\cal C}_X$ 
is the code generated by affine subspaces of codimension 2 in $X$ and
where
$AGL(d,2)$ is the natural subgroup of permutation matrices.  The code
${\cal C}_X$ has parameters $[2^d,1+d+{d \choose 2}, 2^{d-2}]$.  

(ii) Any two labelings as in (i) are conjugate by the action of $H$.

(iii) 
For fixed $q, k$, the sets $A(L,x,q,k)$ (see \ref{layers}) are the elements
of $L[q]$ which are all linear combinations of
$X$ of the form 
$2^{-t} \sum_{x \in A} x\varepsilon_B$, where $A$ is an affine subpace of
$X$ of dimension $a$, $p+a-2t=q$ and $k=p+\lfloor{d\over 2  }\rfloor-2t=
\lfloor{d\over 2  }\rfloor+q-a$ and $\varepsilon_B$ effects sign changes
exactly at indices in $B
\subseteq A$; here $B$ is in the code ${\cal C}_A$, which is spanned by
all $A \cap S$, where $S$ is an affine subspace of codimension 2 in $X$.  

(iv) For a fixed integer $a$, 
the sets $A(L,x,q,k)$ are nonempty exactly
for the indices $k \in I(d,p,q)$  and they are the orbits of
$H$ on
$\mv {L[q]}$. 
\end{thm} 
\pf  
(i): We use notation of \ref{nextbw}.   
We may and do assume that $d\ge 4$.  
There is by induction a basis $X_1$ of $V_1$ contained in $F$ and 
labeling of 
$X_1$ 
by $\Omega_1:= \FF_2^{d-1}$ 
so that we get an
identification of the stabilizer of
$SF(x) \cap V_1$ with 
${\cal E}_{{\cal C}_{\Omega_1}}{:}AGL(d-1,2)$, in
analogous notation.  

The frame is a double basis for the total eigenspace of $E_1$, a maximal
elementary abelian subgroup of a lower group
$R_1$ on $M_1$.  Using our standard diagonal notation
\ref{diagonalnotation}, \ref{bwtclassification}, take involution
$s=t_{12'}$ in dihedral group $D$ and the corresponding subgroup
$E_{12}$ of
$R_{12}$.  Then $s$  interchanges $M_1$ and
$M_2$.  Let $t \in D$ be the SSD
involution associated to $M_1$.  Then $E:=\la E_1 , t \ra$ is a maximal
elementary abelian group in $R$ and its total eigenlattice has the frame
$F$ as a double basis.  Identify $\Omega_1$ with a codimension 1 affine
subspace of $\Omega := \FF_2^d$.  
We define
$\Omega_2$ to be the complement in
$\Omega$ of $\Omega_1$.  Choose any vector $v_0 \in \Omega_2$.  
Let $v_1 \in X_1$ 
be a frame vector labeled by
$0$ and let $v_2:=v_1^s \in X_2:= SF(x) \cap  V_2$.  
Since the action of
$s$ is an isomorphism 
of the transitive $C_H(s)$-sets $X_1$ and $X_2$, the labeling on $X_1$
transfers uniquely to $X_2$ and we translate this labeling to $X_2$
via vector addition by 
$v_0$ to make a labeling of $X_2$ by $\Omega_2$.  The resulting
labeling of
$X$ is uniquely determined (depending on $v_0, s$,
$X_1$).  

From \ref{agld2}, we see that in $\gd d$, a frame stabilizer contains 
a subgroup $J$ isomorphic to $AGL(d,2)$ in the normalizer of $E$ 
which permutes a basis of the eigenlattice.  Its intersection, $K$, 
with a
natural $\gd {d-1}$ subgroup is an analogous $AGL(d-1,2)$ subgroup.  
Let $Z$ be the group generated by $\{ \pm 1_V\}$.   

There are just two $J$-invariant sets of $Z$-orbit representatives in
$F$.  When one of them is restricted to $K$, we get two orbits.  
If $X_1$
is one of these, the other is $X_1^s$ or $-X_1^s$.  
We replace $s$ by $-s$ if necessary to arrange for the other to be
$X_1^s$.  Then $s \in J$.  
The labeling on $X_1$ now extends to all of $X$, 
which is an $H$-invariant
set.

(ii):  Let $\ell, \ell'$ be two labelings for which $H$ is the indicated
monomial group.  We shall transform one to the other by action of $H$. 
Call the {\it domain} of a labeling to be the points of $SF(x)$ which
get a label. 

The stabilizer $H_\ell$ in $H$ of the labeling $\ell$ (equivalently,  of
its domain)  is a complement to the normal
subgroup of sign changes.  Such a subgroup is isomorphic to
$AGL(d,2)$.  
We first note that any two complements are conjugate.  This follows from a
 cohomology argument, \ref{cohomology}.  
From this, we may and do arrange for the two labelings to have the same
domain, which we call $D$.  Since $H$ acts 3-transitively and leaves
invariant a unique Steiner system with parameters $[3,4,2^d]$, addition of
labels of vectors is determined by $H$ once an origin is chosen.  Given an
origin, a partial labeling of $D$ by a basis of $\FF_2^d$
determines the labeling.  Any two such choices lie in one orbit under the
action of
$H$.

(iii) and (iv):  
It  is clear from induction and the form of the types 1, 2 and 3
minimal vectors that a minimal vector has the zoop2 property
\ref{zoop2} with respect to a given sultry frame.  So, the nonempty
sets $A(L,x,q,k)$, for
$k \in I(d,p,q)$, partition
$\mv {L[q]}$.  It remains to 
show that they are orbits for the frame
stabilizer.  

 The action of $AGL(d,2)$ is transitive on affine subspaces
of given dimension.  

  Write $v=v_1+v_2$, where $v_i$ is the
projection to $V_i$, $i=1, 2$.  
Either $v=v_1, v=v_2$ or $v_1\ne 0\ne v_2$ and 
there exist integers $t_i$ and affine subspaces $A_i$ of $X_i$ and
$B_i\in {\cal C}_{A_i}$ so that 
$v_i= 2^{-t_i}\sum_{y \in A_i} y \varepsilon_{B_i}$.
The zoop2 property implies that $t_1=t_2$ and $dim(A_1)=dim(A_2)$. 
Call these common values $t, a$, respectively.  We assume that 
$v_1\ne 0\ne v_2$.  

If there exists an affine hyperplane $X'$ of $X$ so that $U:= supp(v)
\subseteq X'$, we use induction since the $v$ is a minimal vector in the
sBW sublattice of rank $d-1$ supported by $U$.  Suppose that no such $X'$
exists. Then we are in the third case  $v_1\ne 0\ne v_2$ and we use
notation $v_1\ne 0\ne v_2$ as above.  Let $X'$ be any affine hyperplane.  
We claim that
$|U\cap X'|=\half |U|$.  Suppose otherwise.  Then, replacing $X'$ by its
complement, we may assume that $|U\cap X'| <
\half |U|$. Then the sublattice
$S$ of $L$ supported by $X'$ has a vector in $\dual S$ of norm less than
$\half \mu (L)$, a contradiction.  The claim follows.  We get a final
contradiction by using
\ref{rvl}.  
\eop

\begin{rem} \labttr{minvecbe} 
The results \ref{labeling} (iii), (iv),  were proved in \cite{BE}; see 
Th\'eor\`eme I.5, 
Th\'eor\`eme  II.2.  
\end{rem}

\section{Orbits on norm 4 frames in $\leh$. }

We give an application of our theory by giving a short proof that the
Weyl group of $E_8$ has just four orbits on plain frames \ref{plainframe} 
of norm 4 vectors in
$\leh$, equivalently, of
$D_1^8$-sublattices.  This result can be deduced from a classification of
$\ZZ_4$ codes 
\cite{CS}.

\begin{de} \labttr{dinvariant}  Let $L$ be any lattice.  
If $M$ is a sublattice, $2L \le M \le L$, 
the {\it  $d$-invariant
of the frame $F$ (relative to $M$)} is the dimension of the
span of
$F+M/M$.  
Also, we say two plain frames $E, F$ are {\it congruent} if
and only if $E+M=F+M$.

The {\it d-invariant} of a plain frame $F$  
is the dimension of the
subspace of
$L/2L$ spanned by $F+2L$, i.e., the relative $d$-invariant
for $M=2L$.  
\end{de}

\begin{rem}\labttr{din1234}  
Now suppose that $L=\bw 3$.  
The $d$-invariant of a frame 
is a number between 1 and 4 since the image 
is not trivial and spans a
totally singular subspace. 
\end{rem} 

\begin{rem}\labttr{weyle8onrootframes}  It is easy to see that the Weyl
group of $E_8$ is transtive on frames of roots.  This follows from Witt's
theorem since the Weyl group induces the full orthogonal group on
$\leh$ modulo 2 and any frame of roots spans an index 16 sublattice with
all even inner products, hence corresponds mod 2 to a totally isotropic
subspaces with nonsingular vectors.  The next result refers to action of
the proper subgroup $\gd 3$ on frames of roots and norm 4 vectors.  
\end{rem}

\begin{prop}\labtt{G3onnorm2} 
(i)   In the action of $\gd 3$ 
on frames of norm 2 vectors, 
there are four orbits.  
They are distinguished by their $d$-invariants relative to the
sultry twist $L[1]$.

(ii)  In the action of $W_{E_8}$ on frames of norm 4 vectors, 
there are four orbits.  They are distinguished by their $d$-invariants.

\end{prop}
\pf   
(i) It is easy to determine the orbits of $\gd 3$ 
on frames of roots.  
They are represented by the following
vectors with respect to 
$x_1,\dots , x_8$, a standard
orthogonal basis of roots (see \ref{twistede8}):

$F_1: \pm x_1, \dots , \pm x_8$.

$F_2: \pm x_i, i \not  \in A; \half \sum_{j \in A} \pm x_j$, 
where $A$ is a 4-set
of indices representing a Hamming codeword, 
and evenly many signs over $A$
are minus.

$F_3:=   \pm x_i, i  \in B; \half (00aaaa00), \half (0000pqrs), 
\half (00tu00cc)$, 
where $B$ is a
2-set of indices (which we take to be $\{1,2\}$) and 
the indicated partition of the eight indices  into 2-sets has the property
that the union of any two of them is a Hamming codeword.  Also, 
 $a, b, c, p, q, r, s, t, u  \in\{ \pm 1\}$ and where $p=-q, r=-s, t=-u$.

$F_4 := \pm x_1$ and $\pm \half (01111000), \pm \half (0001,-1,110), \pm
\half (0,-1,0,0,1,0,1,1)$.

The proof is an easy 
exercise with the action of the monomial group $H 
\cong 2^7{:}AGL(3,2)$, a subgroup of $\gd 3$, where the group of sign
changes at evenly many indices is indicated by
$2^7$.  Since $\gd 3$ is transitive on roots, an orbit of such  a frame has a
member containing
$x_1$.  We now restrict ourselves to transformations by elements of $H
\le \gd 3$.  If the remaining members of the frame are the
$x_i$, we are in case $F_1$.  If not, one can arrange for the next member
of the frame to be something of the form mentioned in case $F_2$,
supported by a 4-set,
$A$.  If all remaining members of the frame are some $\pm x_j$ or
supported by the same 4-set, we are in the orbit of $F_2$.  If  not,
similar reasoning brings us to case $F_3$ or $F_4$.

One must show that these frames represent different orbits, and
that is accomplished by showing that their images in $L/L[1]$ span
subspaces of dimensions 1, 2, 3 and 4, respectively. (This is
verified by  Smith  canonical  forms, easy to do by hand or
with a software package like Maple): in our notation,
$L[1]$ is the $\ZZ$-span of the $x_i\pm x_j$ and $\half(x_1+\dots +x_8)$. 
These dimensions are the $d$-invariants of the original orbits. 

(ii) 
Let ${\cal O}_i$, for $i=1,\dots ,r$ be the orbits.  
Since $\weh$ induces the full orthogonal group on  
$L/2L$, any orbit has a
representative contained in $L[1]$ since $L[1]/2L$ is a maximal totally
singular subspace.  Now consider the subgroup $\gd 3$, 
which is normalized
by (the nonorthogonal transformation) 
$1-f$, where $f$ is a
fourvolution.  The action of $(1-f)$ takes the set of 240 
roots bijectively to the union of the nonempty sets 
 ${\cal O}_i \cap
L[1]$, and this correspondence preserves orbits of $\gd 3$. We are done by
(i).
\eop

\section{Clean pictures, dirty pictures and transitivity}

We next prove transitivity results for certain kinds of sublattices.   In
particular, we can classify certain scaled embeddings
of $\bw k$ in $\bw d$, for certain $k\le d$.   
See \ref{cleandirty} for the clean and dirty
terminology.

\begin{thm}\labtt{transdirty}  Let $L=\bw d$,  for $d \ge 4$.  There is a
$\gd d$-invariant bijection between sublattices of $L$ which are ssBW of rank
$2^{d-1}$ and noncentral lower involutions, via the SSD correspondence.
\end{thm}
\pf 
Let $M$ be such a sublattice and $t=t_M$ the associated SSD involution.  
Since $t$ normalizes $R$ and has trace 0, it is dirty (see Appendix
A2), whence there is an element $g \in R$ so that $[t,g]=-1$.  We may
arrange for $g$ to be an involution.  Then $g$ interchanges $M$ and
$N:=L\cap M^\perp$, whence
$N$ is a $ss \bw {d-1}$.  
By \ref{dih8onlatt}, 
the condition $det(L^{\pm }(g))=det(M)$ implies that 
$L$ is part of an X-quadruple
$(L,L^+(g),L^-(g),t)$, whence the classification
\ref{bwtclassification} implies that $t$ is  lower. 
\eop

\begin{rem}\labttr{nonssbw} 
There are cases of sublattices $X$ of $\bw d$ of rank $2^{d-1}$ which 
satisfy $L/[X \perp X^\perp]$ elementary abelian, but $X$ is not 
isometric to a scaled $\bw {d-1}$.    For $d=3$, one can take $X$ to be
the sublattice spanned by a root system of type $A_1^4$ which is not
contained in a $D_4$ subsystem.  Such a sublattice is SSD and
corresponds to a SSD involution of trace 0 which is upper with respect to
any conjugate of $\gd 3$ which contains it.  The noncentral involutions of
$R_3$ have trace 0 and fixed point sublattice isometric to $L_{D_4}$.  
\end{rem}

\begin{thm} \labtt{transclean} Suppose that $L=\bw d$ and that $M, M'$ are
sublattices which are the fixed point lattices for clean isometries of
order 2.  If
$rank(M)=rank(M')$, then there is an isometry $g$ of $L$ so that $M'=M^g$.
\end{thm}
\pf 
Such sublattices correspond to SSD involutions with nonzero traces.  Now
use \ref{conjclean}.  
\eop

The following is an application of \ref{transclean}.  

\begin{coro}\labtt{ssbwd-2} 
Suppose that $d\ge 5$ is odd.  Then in $\bw d$ any two ssBW sublattices of
rank
$2^{d-2}$ are in the same orbit under $\gd d$.
\end{coro}
\pf Such sublattices must be SSD.  
\eop

\begin{de} \labttr{ancestral}  
Let $L=\bw d$.  
A {\it first generation sublattice} of $L$ is a sublattice
$L_1$ so that there exists a sublattice $L_2$ 
and an involution $t$ 
so that 
$(L,L_1,L_2,t)\in {\frak X}$.  

A chain of lattices $L=L(0) \ge L(1) \ge \dots \ge L(d)$ is a 
{\it generational chain} if there exists an elementary abelian group $E\le
R$ and a chain of subspaces $E=E(d) > E(1) > \dots > E(0)=\la -1\ra$   
so that for each $k$, $|E(k)|=2^{k+1}$ and $L(k)$ is the total eigenlattice
of
$E(k)$, \ref{eigenlattice}.  

In each $L(k)$, each orthogonally indecomposable summand is a ssBW
sublattice, all of common rank $2^{d-k}$ if $k \le d-2$, and $L(d-1)$ is
a direct sum of isometric rank 1 lattices.  
 Call $L(k)$ a {\it $k^{th}$ generation sublattice} and $E(k)$
{\it  its  defining lower group}.  A sublattice is {\it ancestral} if
it is a 
$k^{th}$ generation sublattice, for some $k$.  
\end{de}

\begin{thm}\labtt{stabancestralfamily} Let $d \ge 4$.  
If $L=\bw d$ and $Z$ is a  $k^{th}$-generation 
sublattice,
$k\le d-2$,  then the stabilizer of $Z$ in $Aut(L)$, is just
$N_{Aut(L)}(E)$, where $E$ is its defining lower  group,
as in 
\ref{ancestral}.
It contains
$\rd d$ and its image in
$\gd d/  \rd d$ is a maximal parabolic which modulo the unipotent
radical has shape 
$GL(k,2)\times
\Omega^+(2(d-k),2)$.  
The $k^{th}$ generation sublattices are in $\gd d$-equivariant bijection
with the elementary abelian subgroups of $\rd d$ which contain $Z(\rd d)$.  
\end{thm}
\pf
The direct summands of $Z$ realize all the linear characters of $E$ which
do not have $-1$  in their kernel.  Thus, $Z$ determines $E$.  By
definition of ancestral sublattices, $E$ determines $Z$.  
\eop

\begin{de}\labttr{ancestorlookalike}
A sublattice of $L=\bw d$ is an {\it $k$-generation ancestor lookalike}
if it is an orthogonal direct sum of $2^k$ copies of ssBW lattices, all of
rank $2^{d-k}$.
\end{de}

The transitivity situation for lookalikes is unclear.  Here is a 
simple result.

\begin{prop} \labtt{1generation} For 
$L=\bw 3$ , there is just one orbit of the automorphism group on third
generation ancestral lookalike 
sublattices and there are four orbits for $\gd 3$.  For $\bw 4$, there are at
least 4 orbits of the automorphism group on third
generation ancestral lookalike 
sublattices.  
\end{prop}
\pf
For the case $L=\bw 3 \cong \leh$, this was covered in \ref{G3onnorm2}. 

Now take the case $L=\bw 4$.  Let $F$ be such a frame.  Then $F+L[1]$
spans a totally singular subspace of $L/L[1]$.  Since $Aut(L)$ induces on
$L/L[1]$ its simple orthogonal group, we may assume that $F$ lies in
the ancestor sublattice $L_1+L_2\cong \rtleh \perp  \rtleh$.    

Since norm 4 elements in $L_1 + L_2$ are indecomposable, we have 
$F=F_1\cup F_2$ where $F_i:=F\cap L_i$.  By using the ideas in the proof
of 
\ref{G3onnorm2}, we find that the
dimension  of the span of $F_2+L_1[1]$ in $L_1/L_1[1]$ 
can be 1, 2, 3 or 4.  
We conclude that the image of $F$ in $L/L[1]$ spans a space of dimension at
most 8 and dimensions 1,2,3 and 4 actually do occur.  This gives a lower
bound of 4 on the number of orbits.  
\eop

\section{The  Ypsilanti lattices }

\def\toteigenlatt#1{L(#1,+)\perp L(#1,-)} 
\def\eigenlattot#1#2{[L_1\perp L_2](#1 , #2 )}
\def\toteigenlattot#1{[L_1\perp L_2](#1 , +)\perp [L_1\perp L_2](#1 , -)}
\def\lot{L_1 \perp L_2}

We now set up a procedure for creating many isometry types of lattices in
sufficiently large dimensions divisible by 8.  Here is a 
rough idea.  We take several  isometric
``good'' lattices (indecomposable, high minimum norm, elementary abelian
discriminant group) and study overlattices $L$ of their 
orthogonal direct sum
$L_1 \perp \dots 
\perp L_s$.  We consider conditions like X (\ref{condx}) but   
without (e).  A suitable concept of avoidance allows us to build
many  lattices $L$ with enough but not too many minimal vectors.  We
gain  enough control over the automorphism groups to get a fairly
high lower bound on the number of isometry types.  

We start with a generalization of the maps $f-1$ where $f$ is a
fourvolution.  

\subsection{Michigan lattices and
Washtenawizations}  

\begin{de} \labttr{specialendo}  A {\it 2-special endomorphism} on a
lattice
$L$ is an endomorphism $p$ so that 

(i) $(xp, yp)=2(x,y)$ for all $x, y \in
L$;  

(ii)   $Lp^2=2L$ (thus, $\half p^2 \in Aut(L)$); 

(iii) there is an integer $r$ so that $\dual L = Lp^{-r}$ ($r$ is called
the {\it duality level}).  

If $L$ has a 2-special endomorphism, call $L$ a {\it 2-special
lattice}.  Call $L$ {\it
normalized} if the duality level is 0 or 1.  
\end{de}  

\begin{rem}\labttr{specialendo2}
A 2-special lattice  is scale-isometric by a
power of a 2-special endomorphism to a normalized lattice.  
\end{rem}

\begin{nota}\labttr{specialnota}   We adapt notations used earlier and
set $L[k]:=Lp^k$, for $k \in \ZZ$.  
When, $C_{Aut(L)}(L[k]/L[k+1])$ is independent of $k \in \ZZ$, we
define
$Lower(L):=C_{Aut(L)}(L/L[1])$ and
$Upper(L):=Stab_{Aut(L)}(L[-1])/Lower(L)$.  
\end{nota}

\begin{nota}\labttr{smv} 
The sublattice of the lattice $L$ spanned by the minimal vectors is
denoted 
$SMV(L)$.  When $L$  has a 2-special endomorphism, define 
$SMV(L,L[1]):=SMV(L)+L[1]/L[1]$ and define 
$mvd(L,L[1])$ to be
the dimension of $SMV(L,L[1])$.  This number is called  {\it the
mv-dimension} and  is positive if
$L \ne 0$.  In case $p$ or $L[1]$ is understood, we write $mvd(L)$ for
$mvd(L,L[1])$ and note that this invariant could depend on choice of
2-special  endomorphism.

Define  the {\it Washtenaw number} or {\it Washtenaw ratio}
of $L \ne 0$ to be the ratio
$$Washtenaw(L):=2\, mvdim(L)/rank(L) = mvdim(L)/dim(L/L[1]) \in
(0,1].$$ 
\end{nota}

\begin{de}\labttr{michigan} A {\it Michigan lattice } is a lattice $M$

(i)  
with a 2-special endomorphism, $p$; 

(ii) $SMV(M)$ has finite index in $M$;    

(iii)  $Aut(M)$ fixes each $Mp^k$, $k \in \ZZ$; 

(iv)  $g\in Aut(M)$ is
trivial on $Mp^k/Mp^{k+1}$ if and only if $g$ is trivial on
$Mp^\ell/Mp^{\ell +1}$, for all $k, \ell \in \ZZ$.  
\end{de}  

Note that a Michigan lattice $L$ is indecomposable if 
$SMV(L)$ is indecomposable.  

\begin{de}\labttr{w2} 
We are given a normalized Michigan lattice $M$ such that
$SMV(M)$ is indecomposable.  
Let $t\ge 3$ be  an integer.

Let $M_1, \dots ,
M_{2^t}$ denote pairwise orthogonal copies of $M$, identified by
isometries 
$\psi_i : M \rightarrow M_i$, 
with 2-special endomorphism $p_i$ corresponding to $p$ by 
 $\psi_i$.  
The direct sum has a 2-special endomorphism, $q$,  which is the direct
sum of the 
$p_i$.

A  {\it  degree $t$ Washtenawization}
of
$M$ is  a lattice $W$  contained in $\QQ \otimes
(M_1\perp
\dots
\perp M_{2^t})$ so that 

(i) $W$ contains $(M_1\perp \dots \perp
M_{2^t})[1-r]$ and is a sublattice of 
$(M_1\perp \dots \perp M_{2^t})[-r]$; ($r$ is the duality level of $M$) and
the quotient 
$M/(M_1\perp \dots \perp M_{2^t})[1-r]$
is elementary abelian of dimension $2^{t-2} \, rank(M)$;   

(ii) For all  $i$, $W \cap (\QQ \otimes M_i)=M_i[1-r]$;   

(iii) $\mu (W)=2^{1-r} \mu (M)$; 

(iv)  $SMV(W)=\sum_{i=1}^{2^t}  SMV(M_i)$ and 
$Washtenaw(W)=\half Washtenaw(M)$; 

(v) $Aut(W)$ has the form $[\prod_{i=1}^{2^t} Lower(M_i)].[Upper(M)
\times Aut({\cal C})]$, where $\cal C$ is an indecomposable
(\ref{defindeccode}) self orthogonal doubly even binary code of length
$2^t$; furthermore,
$Aut(M)$  embeds in
$Aut(W)$ by diagonal action.  

A {\it minimal Washtenawization} is a degree 3 Washtenawization, using
the extended Hamming code (which is essentially the only choice
here).  It is unique up to isometry.  
\end{de} 

\begin{rem} \labttr{dualitywash}   By
\ref{kstlattices},  Washtenawizations are indecomposable, since the
code is indecomposable.  In the notation of
\ref{w2}, the duality level of
$W$ is $1-r$ and $|Upper(W)|$ divides $|Upper(M)|(2^t!)$.  Also,
$Aut(W)$ permutes the set
$\{M_1,
\dots  , M_{2^t}\}$. 
\end{rem}  

\begin{prop}\labtt{w3} For all $t\ge 3$, degree $t$ Washtenawizations
exist. 
\end{prop}
\pf  Let $M$ be a normalized Michigan lattice.   
Take the lattice $W$ between $(M_1\perp \dots \perp
M_{2^t})[1-r]$ and $(M_1\perp \dots \perp
M_{2^t})[-r]$ which corresponds to some indecomposable doubly even self
orthogonal code, $\cal C$ (for example, see \ref{indeccode}).  
Since nonzero code words have weight at least 4, the minimal vectors of
$W$ lie in
$SMV((M_1\perp
\dots
\perp M_{2^t})[1-r])$ (use \ref{projweight}).  

Since $q$ acts diagonally as $p$ on $(M_1\perp \dots \perp
M_{2^t})[1-r]$, the definition of $W$ implies that  the image of 
$SMV((M_1\perp \dots \perp M_{2^t})[1-r])$ in 
$W/Wp$ has dimension $2^{t-1}\,  mvdim(M)$.  
This implies that $Washtenaw(W)=\half Washtenaw(M)$.  

Since
$Aut(W)$ permutes the minimal vectors, it permutes the indecomposable
direct summands of the lattice they generate, which are just  the $2^t$ 
$SMV(M_i)$, which in turn define the  $M_i$ as the summands of $W$
(as abelian groups) which contain the $SMV(M_i)$.  
It follows that $Aut(W)$ is contained in a natural wreath product
$Aut(M_i) \wr Sym_{2^t}$ which permutes $\{M_1,
\dots  , M_{2^t}\}$.  Obviously, $Aut(W)$ contains
a group $G_0$ of the form indicated in \ref{w2}(v).  
Now, use \ref{projweight}(ii) and the fact that $Aut(M)$ leaves each
twist $M[k]$ invariant. 
\eop

\subsection{Overlattices of direct sums of 2-special lattices}

\begin{nota}\labttr{michigan2}  Throughout this section, $M$ is a
normalized 2-special lattice (\ref{specialendo}) and $M_1,
M_2$ are pairwise orthogonal lattices isometric to $M$ with duality level
$r\in \{0,1\}$.  Let $t$ be an isometry of order 2 which interchanges
them.  
\end{nota}  
 
\begin{de}\labttr{admissible}  
The {\it $i^{th}$ admissible component group} $K_i$ is the full general
linear group  on $M_i[-r]/M_i[1-r]$ when  the duality
level of $M$ is 0 and when
the duality level of $M$ is 1, it is the full orthogonal group on the
nonsingular quadratic space
$M_i[-r]/M_i[1-r]$,
$x+M[1-r]\mapsto 2^{r-1} (x,x) (mod \ 2)$. 
\end{de}

\begin{nota}\labttr{fraky} 
Let $d\ge 5$ be an integer and let 
$M_1, M_2$ be isometric normalized 2-special lattices of ranks
$2^{d-1}$ and duality level 1.  Set 
$V_i :=\QQ
\otimes M_i$.     Let ${\frak Y}:= {\frak Y}(M_1[1-r], M_2[1-r] )$ 
denote 
the set of even integral 
lattices $M$ which contain $M_1[1-r]\perp M_2[1-r]$ 
and satisfy  $M \cap V_i =
M_i[1-r]$ for $i=1,2$ and whose projection to $V_i$ is $M_i[-r]$. 
This is a set of rank $2^d$ unimodular lattices.  
(Note differences with
\ref{frakx}, which results in unimodular lattices for ranks $2^d$, $d$
odd only. )

\end{nota}

\begin{rem}\labttr{aboutdiagonallattices} 
A member $L$ of $\frak Y$ is determined 
by an isomorphism of  vector 
spaces $\zeta : M_1[-r]/M_1[1-r] \rightarrow M_2[-r]/M_2[1-r]$, 
namely $L/(M_1[1-r]+M_2[1-r])$ is
just the diagonal in the identification 
of the two $M_i[-r]/M_i[1-r]$ based on $\zeta$. 
We may write $L/(M_1[1-r]+M_2[1-r])=\{ (x+M_1[1-r], (x+M_1[1-r])^\zeta
) | x \in M_1[1-r]\}$.

Conversely, given a linear isomorphism $\zeta$, we get an $L \in {\frak Y}$
by taking the diagonal as above provided (a) when $d-1$ is odd, no
condition; (b) when $d-1$ is even, $\zeta$ is an isometry of nonsingular
quadratic spaces 
$ M_1[-r]/M_1[1-r] \rightarrow M_2[-r]/M_2[1-r]$.  

The reason for the isometry condition in (b) is that the nonsingular
cosets (respectively, the singular cosets) 
of the two $M_i[-r]/M_i[1-r]$ must be
matched to create a diagonal which gives an even lattice $L$.  In (a),
since the two
$M_i[-r]$ are even integral lattices, any matching  by a linear
isomorphism results in an element of $\frak Y$, whence no conditions
are demanded.  The requirement in (b) of taking 
$M_1[-r]/M_1[1-r]$ to $M_2[-r]/M_2[1-r]$ comes from the definition of
$\frak Y$,
\ref{fraky}.  

\end{rem}  

\begin{nota}\labttr{zeta} 
We use
the notations $L \mapsto \zeta (L), \zeta \mapsto L(\zeta)$ to express
the bijection between $\frak Y$ and  such isomorphisms.  

Such $\zeta$ are in bijection with $K_1$ and with $K_2$ (see
\ref{admissible}) by
$\zeta
\mapsto \zeta_i \in K_i$, where the latter are defined by the formulas
 $ \zeta : x+M_1[1-r] \mapsto (x^t+M_2[1-r])^{\zeta_2} =
y^t+M_2$, where $y +M_2=(x+M_1)^{\zeta_1 t} $, where $t$ is as in
\ref{michigan2}.  Call
$\zeta_i$ {\it the
 $K_i$-component of $\zeta$}, or of $L=L(\zeta )$.  
\end{nota}

\subsection{Avoidance}  

\begin{de}\labttr{avoid}  We say that two subspaces of a vector space
{\it avoid} each other if their intersection is 0.  If $g:V\rightarrow V'$
is an invertible linear transformation, $W\le V$ and $W'\le V'$, we
say that
$g$ is a
$(W, W')$-{\it avoiding map}  if $W^g\cap W'=0$.  Let 
 $A(W_1,W_2)$ be the set of avoiding maps.  
\end{de}

We need some terminology for discussing asymptotic behavior.

\begin{nota}\labttr{dtl}    Suppose that $f(x)$ is a
real-valued function on $(0,\infty )$.  The  {\it dominant
term} in $f(x)$ (abbreviated
$DT(f(x))$ is 
the expression of the form $a_0log_2(x)^{a_1} 2^{a_2x} x^{a_3}$ which is
asymptotic to
$f(x)$ (the $a_i$ are constants).  We may indicate dependence on the
variable
$x$ by
$DT_x$.  (This definition applies to a limited family of real-valued
functions, but suffices for our purposes.) 

Similarly, if $f$ is as above, we define the {\it dominant term of the
logarithm} ($DTL$ or $DTL_x$) of $2^{f(x)}$ to be $DT(f(x))$.  
For example, $$DTL(2^{(0.43)log_2(2x-3) 
2^{3x-4}+2^{2x}-log_2(x+1)^5 x^3-log_2(x)^7(x^2+1)})=
\hbox{${0.43 \over 16}$}   \
log_2(x) 2^{3x-4}.$$  
\end{nota}

\begin{prop}\labtt{avoidercount}  Suppose that $a \le b$ are positive
integers.  Suppose that $V:=\FF_2^{2b}$ has a maximal Witt index
nonsingular  quadratic form and that $W_1$ and $W_2$ are two
$a$-dimensional totally singular subspaces.  We set $q:={a\over b}$ and
think  of $q$ as a constant and $a$  as a function of
$b$.  

(i) Let $H$ be the stabilizer in $O(V)$ of $W_1$.  Then, 
$$DTL_b(|H|)=
DT_b({\hbox{$\half$} a(3a-1)+2(b-a)b})=b^2(2-2q+\hbox{${3\over
2}$}q^2).$$ 

(ii) For an integer $k$, let  $A(W_1,W_2;k)$ be the set of avoiding maps
as in \ref{avoid} so that $dim(W_1^g\cap W_2^\perp )=k$.  Then 
$A(W_1,W_2;k)$  is nonempty precisely for $k=0, 1, \dots , min\{a,b-a\}$
and for each such $k$, $A(W_1,W_2;k)$ is a regular orbit for the action of
$H$.  
\end{prop} 
\pf  (i): We have $DTL_k(|\Omega^+(2k,2)|)={2k^2-k}$.  We may assume
$W_1=W_2$. Let
$H$ be the subgroup of the orthogonal group which fixes $W_1$ globally.
It follows from \ref{aboutomega} that $DTL(|H|)$ is the DTL of $\half
a(3a-1)+2(b-a)b-(b-a)$.

(ii) : By Witt's theorem, 
two nonavoiding maps $g, g'$ are in the same $H$-orbit if the dimensions of
the images of $W_1$ under $g, g'$ intersect $W_1^\perp$ in spaces of the
same dimensions.  All $H$-orbits are regular.  For a nonempty
$A(W_1,W_2;k)$, we have $k \le dim(W_1)=a$ and since the image of an
avoiding $W_1$ in $V/W_2^\perp$ has dimension at most
$a=dim(V/W_2^\perp)$, we have
$k+a=k+dim(W_2)
\le b$, the dimension of any maximal totally isotropic subspace.  The
value $k=min\{a,b-a\}$ can be achieved.  
\eop 

\begin{coro} \labtt{q}  We use the notations of \ref{avoidercount} and
assume that $q\le \half$.   Then $$DTL_b(|A(W_1,W_2)|)=
\upsilon (q) log_2(b)b^2, \hbox{where } \upsilon (q):=(2-2q+\hbox{${3\over 2}$}
q^2).$$
\end{coro}  
\pf  Note that $q\le \half$ means $a=min\{a, b-a\}$ in \ref{avoidercount}.
\eop

\subsection{Down Washtenaw Avenue to Ypsilanti} 

We next create large families of lattices in  dimensions $2^d >> 0$.  

\begin{de} \labttr{ypsicousins}  
Let $W$ be a normalized Michigan lattice which has duality level $r=1$
and Washtenaw ratio $q\le \half$.  

Take orthogonal copies $M_1, M_2$ of
$W$ and consider the set ${\frak Y}:={\frak Y}(M_1, M_2)$ as in
\ref{fraky}.  Consider the associated maps $\zeta (L), L \in  {\frak Y}$
(see \ref{zeta}) which are avoiding maps \ref{avoid} for the subspaces
$SMV(M_1[-1],M_1), SMV(M_2[-1],M_2)$ of $M_1[-1]/M_1$, $M_2[-1]/M_2$,
respectively.  The corresponding lattices form a subset 
${\frak Y}_{av}(M_1,M_2)$ of ${\frak Y}(M_1,M_2)$ in the notation of
\ref{fraky}.   Their
ranks are $2\, rank(W)$.  They are called {\it Ypsilanti lattices}.   Let 
$IsomTypes(M_1,M_2)$ be the set of isometry types of lattices in ${\frak
Y}_{av}(M_1,M_2)$.  

When $W$ is a Washtenawization of a BW lattice, the Ypsilanti lattices  
of rank $2^d=2\, rank(W)$ are
called the {\it Ypsilanti  cousins} of $\bw d$. 
\end{de}

\begin{lem}\labtt{smvypsi} We use the notations of \ref{ypsicousins}.  

(i) If $N \in {\frak Y}_{av}(M_1,M_2)$, $SMV(N)=SMV(M_1)\perp SMV(M_2)$.

(ii) $N \in {\frak Y}_{av}(M_1,M_2)$ is indecomposable.

\end{lem}
\pf (i) 
Obviously, $\mu (N) \ge \mu  (M_i), i=1,2$.  Consider a vector
$x=x_1+x_2\in N \setminus (M_1\perp M_2)$.  Then the $x_i$ have norms at
least $\mu (M_i[-1])=\half \mu (M_i)$.  For $(x,x)$ to equal $\mu (M_i)$,
we need  $x_i$ to be a minimal vector of $M_i[-1]$ for $i=1, 2$.   This
is not the case since $N$ was defined with an avoiding map.

(ii) Use \ref{indecoversmv}. 
\eop

\begin{lem}\labtt{seriesq}  Suppose that we are given $q=2^{-j}$ for some
$j > 0$.  For all $k\ge 5+3j$, there exists a Michigan lattice $W(k)$ so
that 
$rank(W(k))=2^k$, $Washtenaw(W(k))=q$ and the duality level of $W(k)$
is 1.   We may also arrange for 
 $\mu (W(k)) = 2^{ 1-r+\lfloor {j\over 2}\rfloor + \lfloor {e\over
2}\rfloor}$, where $e=k-3j$ if $k$ is even and   $e=k-3j-1$ if $k$
is odd.   
\end{lem}
\pf
If we start with $\bw e$ and perform the minimal Washtenawization
procedure
$s$ times, we get a lattice $W(e,s)$ of rank $2^{e+3s}$.  
We may take a degree 4 Washtenawization to $W(e,s)$ and get
a lattice $W'(e,s+1)$ of rank $2^{e+3s+4}$.  

Each Washtenawization changes duality level.  
We define $W(k)$ according to the following cases. When $k$ is even, we
require $k-3j\ge 4$, which means $k$ is at least 8.  When $k$ is odd, we
require $k-3j-1\ge 4$, which means that $k$ is at least 9.  

$W(k):=\begin{cases} 
W(k-3j,j) & k \hbox{ even}; \cr 
W'(k-3j-1,j) & k \hbox{ odd}.   \cr
\end{cases}
$
\eop

\begin{de}  \labttr{washtenawseries} We call a sequence  of lattices as
in \ref{seriesq} {\it the $j$-Washtenaw series, for the fixed ratio
$q=2^{-j}$}.   It starts at rank $2^{5+3j}$.  The isometry types of
certain members of the series depend on choice of indecomposable doubly
even code of length 16.  Ypsilanti cousins associated to such series are
called {\it Ypsilanti $j$-cousins}.  The set of such isometry types is
denoted
$Ypsi(2^d,j)$.  
\end{de}

\begin{lem} \labtt{propertiesypsi}  We use the notations of
\ref{ypsicousins}, \ref{seriesq} and let $W(k)$ be the Washtenaw series.  

(i) 
If $N$ and $N'$ are
two cousins of rank $2^d$, $Isom(N,N')$ is contained in the group
$G_0(e,h)$ of orthogonal transformations which stabilize $L_1 \perp
\dots \perp L_{2^h}$,  the indecomposable direct summands of
$SMV(N)=SMV(N')$ (in fact, the 
$L_i$ are the pairwise isometric scaled Barnes-Wall lattices, of rank
$2^e$, on which the Washtenawizations $M_1, M_2$ were based; the
notation means $d=k+1=e+h$, wih $h=3j$ or $3j+1$).  

(ii) $DTL_d(|G_0(e,h)|)$  is bounded above by a constant
times
$d^2$.  
\end{lem} 
\pf (i) 
Given $N, N' \in {\frak Y}_{av}(M_1,M_2)$, an
isometry of
$N$ to
$N'$ takes $SMV(N)$ to $SMV(N')$.  Both of these equal 
$SMV(M_1\perp M_2)$.

(ii)  This follows from $DTL_f(|\Omega^+ (2f,2)|)=2f^2$ and boundedness
of $h$.  
\eop

\begin{nota}\labttr{upsilon} When $W(k)$ runs through 
the $j$-Washtenaw series \ref{washtenawseries}, we let $\Upsilon
(2^d,j):= |Ypsi(2^d,j)|$.  
\end{nota}

\begin{lem}\labtt{upsilon2} $\Upsilon (2^d,j)\ge |{\frak
Y}_{av}(M_1,M_2)|/|G_0(e,h)|$, whence 

\noindent $DTL_d(\Upsilon(2^d,j)) \ge   {1\over 16}  \upsilon
(2^{-j})d\, 2^{2d}$,  as in  \ref{q}.  
\end{lem}   
\pf   In \ref{q}, take $b=2^{d-2}$, because the admissible
component group \ref{admissible} is $O^+(2^{d-1},2)$ since the
duality level has been arranged to be 1.  Then use 
\ref{propertiesypsi}(ii)  and 
\ref{upsilon}.\eop

\begin{rem}\labttr{upsilonallq}  For a fixed large value of $d$, we
can make the families $Ypsi(2^d,j)$, $q=2^{-j}$, for all 
$1\le j \le \lfloor {d-5\over 3}\rfloor$.  This would make roughly
$d/3$ times as many as one of the $Ypsi(2^d,j)$, so would not
increase the DTL.
\end{rem}  

We summarize our counting in dimensions $2^d$.  

\begin{thm}\labtt{ypsicount}  For any  $j>0$, the number of Ypsilanti
lattices in dimension $2^d$ has DTL at least $\sixteenth \upsilon
(2^{-j})  d\, 2^{2d}$.  In particular, the number of indecomposable
even unimodular lattices in dimensions $2^d$ has DTL at least $c\  d\,
2^{2d}$, for any $c \in (0,\eighth )$.  
\end{thm}

\begin{rem} \labttr{smalld}   With a bit more work, we could define
lattices like Ypsilanti cousins for $d < 9$, though 
we would not expect them to represent more than a fraction of 
$mass(2^d)$ isometry types.  In dimension 32, the mass formula
gives value about
$10^7$ and the number of isometry types (still not known) 
has been
bounded below by about
$10^{10}$ (see \cite{K}).  
\end{rem}

\subsection{From dimensions $2^d$ to arbitrary dimensions}

\begin{nota} \labttr{ypsin}  For an integer $n>0$ divisible  by 8,
let $2^d$ be the largest power of 2 less than or equal to $n$.  
Fix some $q^{-j}$, $j > 0$.  
Let  
$Ypsi(n,j)$ be the set of isometry types of even integral unimodular
lattices which contain a Ypsilanti $j$-cousin of rank $2^d$ as an
orthogonal direct summand.  Clearly, $\Upsilon (n,j) := |Ypsi(n,j)|
\ge
\Upsilon(2^d,j)$.  
\end{nota}

\begin{coro} \labtt{arbdim} We use the notation of \ref{ypsin}.  For any
constant $c \in [0, {1 \over 32 }  )$, we take $j>0$ so that
$q=2^{-j}$ satisfies $2-q+{3 \over 2}q^2 > 64 c$.  

Then 
$log_2(\Upsilon (n,j)|)\ge   c\, log_2(n)
\, n^2$.  
\end{coro}  
\pf 
Take the integer $d$ which satisfies $2^d \le n < 2^{d+1}$.  Then $d
< log_2(n) \le d+1$ and $2^d > {n \over 2}$.  We have 
$\Upsilon (n,j)| \ge \Upsilon(2^d,j)$ and 
$DTL (\Upsilon (n,j)|)\ge DTL(\Upsilon(2^d,j) ) >  4c\, d\, 2^d\ge
{4c}(log_2 (n)-1)({n\over 2})^2$ whose DT is at least ${c} \
log_2(n)n^2$.  
\eop

\subsection{Number  of Ypsilanti cousins compared with the mass formula} 

\begin{nota}\labttr{massnota} 
We follow the notations of \cite{Se}, pp. 54, 90, except we write
$mass(n)$ instead of ``$M_n$''.  Stirling's
formula ($n!
\sim  n^{n+\half }e^{-n} (2\pi)^\half$) implies that
$DTL(n!)=log_2(n)n$.  
Let $B_j$ be  the  $j^{th}$ Bernoulli number.  
Let $n\in  8\ZZ$, $k:={n\over 8}$. 
\end{nota}

\begin{prop}\labtt{dtlmass} $DTL_n(mass(n))= \fourth log_2(n)\, n^2$. 
\end{prop} 
\pf
We have 
$mass(n)={B_{2k}\over 8k} \prod_{j=1}^{4k-1} {B_j\over 4j}$.  
Because $\zeta(2j)\in
(1,2)$ for all $j\ge 1$, the formula
$B_j=2\zeta(2j).(2j)!/(2 \pi)^{2j}$ shows that 
$DTL_j(B_j)=DT(log_2(2j) 2j)$.  

We have 
$$log_2(mass(n))=log_2(B_{2k})+\sum_{j=1}^{4k-1} log_2(B_j) 
-(8k+1) - log_2((4k-1)!)-log_2(k).$$ 
Since $DTL_k(B_{2k})=DTL_k((4k)!)$, $DTL_n(mass(n))=
DT_n( \sum_{j=1}^{4k-1}   log_2(2j) 2j )$.   
The latter summation can be thought of as   Riemann  sums,   which can be
estimated with integrals (think of 
$\int 2x\, ln(2x)\, dx = 
\int 2x\, ln(x)\, dx + ln(2) \int 2x\, dx  = 
x^2\, ln(x) -\half x^2 + ln(2) x^2 + c$, which 
has dominant term $x^2\, ln(x) $).    We conclude that
$DTL_n(mass(n))=DT_n(log_2(8k) (4k)^2 ) = \fourth   log_2(n)n^2$.  
\eop 

\begin{prop}\labtt{fsgrat} 
For positive integers $n, q$, define $A(n,q):=\sum_{i \ge 0} \lfloor
{n \over q^i (q-1) } \rfloor$ and let $\cal P$ be the set of prime
numbers at most $n+1$.  Set $f(n) := \prod_{q \in {\cal P}}
q^{A(n,q)}$.  Then a finite subgroup of $GL(n,\QQ )$ has order
dividing
$f(n)$.
\end{prop}
\pf  
This is a result of Minkowski \cite{Mink}.  See the discussions in
exercises for Section 7 of 
\cite{Bour}.  
\eop 

\begin{lem}\labtt{dtlmink} $DTL_n (f(n)) = n\, log_2(n)$.  
\end{lem}
\pf
Well known?  A proof may be deduced from   \cite{NZM}, Th. 8.8(b),
p. 369.   
\eop

\begin{rem} \labttr{nuun}
The DTL $n \, log_2 (n)$ is small compared to $DTL(mass(n))$.   It
follows that the DTL of the number of isometry types of rank $n$ even
unimodular lattices is the same as that of $DTL(mass(n))$.  
\end{rem}

We summarize:  

\begin{coro}\labtt{summaryypsi}  For any $a \in (0,\eighth) $, 
there is an integer $j$ so that 
$DTL_n ( \Upsilon (n,j))\ge a\cdot  DTL(mass(n))$.  
Furthermore, when $n$ is a power of 2, 
and 
$b \in (0,\half ) $, 
there is an integer $j$ so that 
$DTL_n ( \Upsilon (n,j))\ge b\cdot  DTL(mass(n))$. 
\end{coro}

\begin{rem}\labttr{newlowerbound}   
We conclude  with some numerical  comparisions.  
\bigskip

\vfill  \eject

\centerline{\bf  Asymptotics for $\Upsilon(2^d,j)$, $\Upsilon (n,j)|$ and
$mass(n)$.  }  
$$\begin{matrix}  j & q=2^{-j}&\upsilon (q) = \hbox{constant
coefficient of } &  
\hbox{Lower bound for}   \cr 
&&16\, DTL(\Upsilon(2^d,j) ) (\hbox{see
\ref{upsilon2}}) &  DTL(\Upsilon (2^d,j)|)/DTL(mass(2^d)) \cr
&&& \cr 
           1&     .5000000000 & 1.375000000 & .3437500000\cr
           2&     .2500000000 & 1.593750000 & .3984375000\cr
            3&    .1250000000 & 1.773437500 & .4433593750\cr
           4&     .06250000000 & 1.880859375 & .4702148438\cr
           5&     .03125000000 & 1.938964844 & .4847412109\cr
            6&    .01562500000 & 1.969116211 & .4922790527\cr
           7&    .007812500000 & 1.984466553 & .4961166382\cr
            8&   .003906250000 & 1.992210388 & .4980525970\cr
            9&   .001953125000 & 1.996099472 & .4990248680\cr
           10&    .0009765625000 & 1.998048306 & .4995120764\cr
\end{matrix}
$$ 
\end{rem}

\section{Appendices}

\subsection{A1.  Group orders}

\begin{prop}\labtt{aboutomega} (i) 
For $q$ a power of 2, the order of $\Omega^+(2n,q)$ is 
$q^{n(n-1)} (q^n-1) \prod_{i=1}^{n-1} (q^{2i}-1)$.  
 
(ii) 
The stabilizer in $\Omega^+(2n,q)$ of an isotropic point has shape 
$q^{2(n-1)}{:}[\Omega^+(2(n-1),q) \times q-1]$.

(iii)  
The stabilizer  in $\Omega^+(2n,q)$ of a maximal totally singular
subspace has shape
$q^{{n \choose 2}}{:}GL(n,q)$, and this is a maximal subgroup.  

(iii)  The stabilizer in $\Omega^+(2n,q)$ of a totally singular
subspace of dimension $m <  n$ has the form $RL$, where the unipotent
radical has order $2^{{m
\choose 2}+2m(n-m)}$ and 
$L\cong \Omega^+(2(n-m),2)\times GL(m,2)$.  These are maximal
subgroups.
\end{prop}  

\pf  These are well-known properties of the orthogonal groups.  Proofs 
may be obtained from \cite{Car, GrElAb}.
\eop

\subsection{A2.  $Aut^0(\exsp d \varepsilon)$, $Out^0(\exsp d
\varepsilon)$ and $BRW^0(2^d,\varepsilon)$.  }

Basic theory of extraspecial groups extended upwards by their
outer automorphism group has been
developed in several places.  We shall use 
\cite{Grex, GrMont, GrDemp, GrNW, Hup, BRW1, BRW2, B}.

\begin{nota}\labttr{brw}  Let $R\cong 2^{1+2d}_\varepsilon$ be an
extraspecial group  which is a subgroup of $G:=GL(2^d,\FF)$, for a field
$\FF$ of characteristic 0. Let $N:=N_G(R)\cong \FF^\times.2^{2d}
O^\varepsilon (2d,2)$.  The {\it Bolt-Room-Wall group} is a subgroup of
this of the form $2^{1+2d}_\varepsilon.\Omega^\varepsilon (2d,2)$.  
If $d\ge 3$ or $d=2,
\varepsilon =-$,  $N'$ has this property.  For
the excluded parameters, we take a suitable subgroup of such a group for
larger $d$.  
We denote this group by $BRW^0(2^d,+)$ or $\dg d$.  It is uniquely
determined up to conjugacy in  $G$ 
by its isomorphism type if  $d\ge 3$ or $d=2,
\varepsilon =-$.  It is conjugate to a subgroup of $GL(2^d,\QQ)$ if
$\varepsilon = +$.  Let $R=\rd d$ denote $O_2( \gd d )$.  We call
$R_d$ {\it the lower group} of $BRW^0(2^d,+)$ and call $G_d/R_d$ {\it
the upper group} of $BRW^0(2^d,+)$.

For $g\in N$, define 
$C_{R \ mod \ R'}(g):=\{ x \in R | [x,g]\in R'\}$, 
$B(g):=Z(C_{R \ mod \ R'}(g))$ and let $A(g)$ be some subgroup of 
$C_{R \ mod \ R'}(g)$ which contains $R'$ and complements $B(g)$ modulo
$R'$, i.e., 
$C_{R \ mod \ R'}(g)=A(g)B(g)$ and $A(g)\cap B(g)=R'$.  
Thus, $A(g)$ is extraspecial or cyclic of order 2.  
Define 
$c(d):=dim(C_{R/R'}(g))$, $a(g):=\half |A(g)/R'|$, $b(g):=\half
|B(g)/R'|$. 
 Then $c(d)=2a(d)+2b(d)$.  
\end{nota}

\begin{coro}\labttr{agld2}  Let $L$ be any $\ZZ$-lattice invariant under
$H:=BRW^0(2^d,+)$.   Then $H$ contains a subgroup $K\cong AGL(d,2)$
and $L$ has a linearly independent set of vectors $\{ x_i | i \in
\Omega\}$ so that there exists and identification of $\Omega$ with
$\FF_2^d$ which makes the $\ZZ$-span of 
$\{ x_i| i \in\Omega \}$ a  permutation module for $AGL(d,2)$ on
$\Omega$.  
\end{coro}
\pf
In $H$, let $E, F$ be  maximal elementary abelian subgroups and let $K$
be their common normalizer.  It satisfies $K/R \cong GL(d,2)$.  Now, let
$z$ generate $Z(R)$ and let $E_1$ complement $\la z \ra$ in $E$ and 
$F_1$
complement $\la z \ra$ in $F$.  The action of $K$ on the hyperplanes of
$E$ which  complement $Z(R)$ satisfies $N_K(E_1)F=K,  N_F(E_1)=Z(R)$. 
Now consider the action of $N_K(E_1)$ on the hyperplanes of $F$ which
complement
$Z(R)$.  We have that $K_1:=N_K(E_1)\cap N_K(F_1)$ covers
$N_K(E_1)/E$.  Therefore, $K_1/Z(R)\cong GL(d,2)$.  Let $K_0$ be the
subgroup of index 2 which acts trivially on the fixed points on $L$ of
$E_1$, a rank 1 lattice.  So, $K_0\cong GL(d,2)$.    Let $x$ be a basis
element of this fixed point lattice.  Then the semidirect product 
$F_1{:}K_0$ is isomorphic to $AGL(d,2)$ and $\{ x^g|g \in F_1\}$ is a
permutation basis of 
its $\ZZ$-span.  
\eop

\begin{de}\labttr{cleandirty}  We use the notation of \ref{brw}.   An
element
$x \in N$ is  {\it dirty} if there exists $g$ so that 
$[x,g]=xz$, where $z$ is an element of  order 2 in the center.
If $g$ can be chosen to be of order 2, call $x$ {\it really dirty} or
{\it extra dirty}.  If $x$ is not dirty, call $x$ {\it clean}.
\end{de}

\begin{lem}\labtt{outer}  Let $\FF_2^{2d}$ be equipped with a
nondegenerate quadratic form with  maximal Witt index.  The set of maximal
totally singular subspaces has two orbits under $\Omega^+(2d,2)$ and these
are interchanged by the elements of $O^+(2d,2)$ outside 
 $\Omega^+(2d,2)$.  
\end{lem}
\pf
This is surely well known.  For a proof, see \cite{GrElAb}. 
\eop 

\begin{thm}\labtt{traces}
We use the notation of \ref{brw}, \ref{cleandirty}.  Let $g \in N$. 
Then $Tr(g)=0$ if and only
if $g$ is dirty.  Assume now that $g$ is clean and has finite order. 
Then 
$Tr(g)=\pm 2^{a(g)+b(g)}\eta$, where $\eta$ is a root of unity.   
If $g\in BRW(d,+)$,
we may take $\eta = 1$. 
Furthermore, every coset of $R$ in $BRW(d,\varepsilon )$ contains a
clean element and if $g$ is clean, 
the set of clean elements in $Rg$ is 
just $g^R\cup -g^R$.  
\end{thm} 
\pf 
\cite{GrMont}.  
\eop

\begin{lem} \labtt{conjorthoginvols} Suppose that $t, u$ are involutions
in $\Omega^+(2d,2)$, for $d \ge 2$.  Suppose that their commutators on the
natural module $W:=\FF_2^{2d}$ are totally singular subspaces of the same
dimension, $e$.  Suppose that $e<d$ or that $e=d$ and that $[W,t]$ and
$[W,u]$ are in the same orbit under $\Omega^+(2d,2)$.   Then $t$ and $u$
are conjugate.  
\end{lem}
\pf 
Induction on $d$.
\eop

\begin{coro}\labtt{conjclean}  Suppose that $t, u$ are clean
involutions in
$H$ with $Tr(t)=Tr(u)\ne 0$.  Then $t$ and $u$ are conjugate in $\gd d$.
\end{coro}
\pf We may assume that $t,u$ are noncentral.  
These involutions are not lower and have the same dimension of fixed
points on $R/R'\cong \FF_2^{2d}$.  Let $T, U\le R$ be their respective
centralizers in $R$.  Since both $t, u$ are clean, $[R,t]$ and $[R,u]$
are elementary abelian subgroups of $T, U$, respectively.  
From \ref{conjorthoginvols}, we deduce that $Rt$ and $Ru$ are conjugate
in $\gd d$.  We may assume that $Rt=Ru$.  Now use \ref{traces} to deduce
that $t$ is $R$-conjugate to $u$ or $-u$. The trace condition implies
that $t$ is conjugate to $u$.  
\eop 

\begin{rem}\labttr {splitting}  The extension 
$1\rightarrow \rd d \rightarrow \gd d \rightarrow \Omega^+(2d,2) \rightarrow
1$ 
is nonsplit for $d\ge 4$.  This was proved first in \cite{BRW2}, then
later in \cite{BE} and in \cite{Grex} (for both kinds of
extraspecial groups, though with an error for
$d=3$; see
\cite{GrDemp} for a correction).  The article \cite{Grex} gives a
sufficient condition for a subextension 
$1\rightarrow \rd d \rightarrow H \rightarrow H/\rd d \rightarrow
1$ to be split, and there are interesting applications, e.g. to the
centralizer of a 2-central involution in the Monster.  A general
discussion of exceptional cohomology in simple group theory is in
\cite{GrNW}.  
\end{rem}

\subsection{A3. Indecomposable integral representations for a group of
order 2}  
\begin{prop}\labttr {indec2}  Let $G$ be a cyclic group of order $2$ and
$M$ a finitely generated $\ZZ$-free $G$-module.  Then $M$ is a direct sum
of modules isomorphic to $\ZZ[G]$, the group algebra; the $\ZZ$-rank  1
trivial module; the $\ZZ$-rank 1 nontrivial $G$-module.
\end{prop}
\pf  \cite{CR}, Section 74.  The case where $G$ has order any prime
number is treated.  \eop



\begin{thebibliography}{GRC99}




\bibitem[1]{BW} E. S. Barnes  and G. E. Wall, Some extreme forms
defined in terms of abelian groups, JAMS 1 (1959), 47-63.   

\bibitem[2]{BRW1} Beverly Bolt, T. G. Room and G. E. Wall, 
On the Clifford Collineations, Transform and Similarity Groups, I. 
Journal of the Australian Mathematical Society, 2, 1961, 60-79.  


\bibitem[3]{BRW2} Beverly Bolt, T. G. Room and G. E. Wall, 
On the Clifford Collineations, Transform and Similarity Groups, II.
Journal of the Australian Mathematical Society, 1961, 80-96.  

\bibitem[4]{B} Beverly Bolt, T. G. Room and G. E. Wall, 
On the Clifford Collineations, Transform and Similarity Groups, III;
Generators and Relations,  Journal of the Australian Mathematical Society,
1961

\bibitem[5]{Bour} N. Bourbaki, \'Elements de Math\' ematique,
Groupes  et alg\`ebres de Lie, Chapitres 2 et 3, Diffusion C.C. L.
S., Paris, 1972.  


\bibitem[6]{BE} Michel Brou\'e and Michel Enguehard, Une famille
infinie de formes quadratiques enti\`ere; leurs groupes d'automorphismes,
Ann. scient. \'Ec. Norm. Sup., $4^{eme}$ s\'erie, t. 6, 1973, 17-52. 

\bibitem[7]{Car}  Roger Carter, Simple Groups of Lie Type,
Wiley-Interscience, London (1972).  


\bibitem[8]{CS} John Conway and Neil Sloane, Self-dual codes over the integers
modulo 4, Journal of Combinatorial Theory, Series A 62, 30-45 (1993).  


\bibitem [9]{CR} Charles Curtis and Irving Reiner, Representation Theory 
of Groups and Associative Algebras, Interscience, 1962.  


\bibitem[10]{DGH} C. Dong, R. Griess. Jr. and G. Hoehn, Framed vertex operator algebras, codes and the moonshine module,
{\em Comm. Math. Phys.} {\bf 193} (1998), 407-448.


\bibitem[11]{Gor} 
 Daniel Gorenstein, Finite Groups, Harper and Row, New York, 1968.  



\bibitem[12]{GrDemp} Robert L. Griess, Jr., On a subgroup of order 
$2^{15}|GL(5,2)|$  in $E_8 (C)$,  the Dempwolff
group and $Aut(D_8 \circ D_8 \circ D_8)$ , J. Algebra, 40, 1976, 271-279.


\bibitem[13]{GrElAb}   Robert L. Griess, Jr., Elementary abelian
subgroups of algebraic groups, Geometria Dedicata, {\bf 39}, 253-305, 1991.  


\bibitem[14]{GrE8} Robert L. Griess, Jr.,   Positive definite
lattices of rank at most 8,  Journal of Number Theory, 103
(2003), 77-84.  


\bibitem[15]{Grex}Robert L. Griess, Jr.,  Automorphisms of extra
special groups and nonvanishing degree 2 cohomology, Pacific J. Math., 48,
403-422, 1973.

\bibitem[16]{GrNW}  
 Robert L. Griess, Jr., 
Sporadic groups, code loops and nonvanishing cohomology, J. Pure Appl.
Algebra, 44, 1987, 191-214.



\bibitem[17]{GrMont}    Robert L. Griess, Jr., 
The monster and its nonassociative algebra, in Proceedings of the Montreal
Conference on Finite Groups, Contemporary Mathematics, 45, 121-157, 1985, 
American Mathematical Society, Providence, RI.



\bibitem[18]{G12} Robert L. Griess, Jr.,   Twelve Sporadic Groups,
Springer Verlag, 1998.  


\bibitem[19]{POE} Robert L. Griess, Jr, Pieces of Eight, Advances in
Mathematics, 148, 75-104 (1999).


\bibitem[20]{Ham} R. W.  Hamming, Error detecting and error correcting
codes, Bell  Syst. Tech. J. 29 (1950), 147-160.  MR 12:35.


\bibitem[21]{Hup} Bertram Huppert,  Endliche Gruppen I, Springer Verlag,  
Berlin, 1968. 


\bibitem[22]{K} Oliver King, 
A mass formula for unimodular lattices with no roots. (English. English
summary)
Math. Comp. 72 (2003), no. 242, 839--863 (electronic).
11H55 (11E41) MR1954971 (Review)

\bibitem[23]{Kneser} M. Kneser, Theorie der Kristalgitter, Math. Ann.
127,  105-106 (1954).  


\bibitem[24]{MS}
Jesse MacWilliams and Neal Sloane, The Theory of Error Correcting Codes,
North-Holland, 1977.

\bibitem[25]{MH} Milnor and Husemoller,  Symmetric Bilinear Forms,
Ergebnisse der Mathematick und Ihrer Grenzgebiete, Band 73, Springer
Verlag, New York, 1973.  

\bibitem[26]{Mink} H. Minkowski, 
 Zur Theorie der positiven quadratischen Formen,
Jour. f\"ur die  reine und angew. Math.  101 (1887), 196-202.  

\bibitem[27]{NZM} Ivan Niven, Herbert S. Zuckerman, Hugh L. 
Montgomery, An Introduction to the Theory of Numbers, Fifth Edition,
Wiley, New York, 1991.  


\bibitem[28]{RVL} B. L. Rothschild and J. H. van Lint, Characterizing
finite subspaces, J. comb Theory, 16A, (1974) 97-110.  

\bibitem[29]{Se} Jean-Pierre Serre, A Course in Arithmetic, Springer
Verlag, Graduate Texts in Mathematics 7, 1973.  



\end{thebibliography}
\end{document}